\documentclass[nonblindrev]{informs3}

\usepackage{natbib}
 \bibpunct[, ]{(}{)}{,}{a}{}{,}%
 \def\bibfont{\small}%

\usepackage{subcaption}
\usepackage{color}
\usepackage{multicol}
\usepackage{booktabs}

\TheoremsNumberedThrough     

\EquationsNumberedThrough    

\begin{document}


\RUNAUTHOR{Han, He, and Oh}

\RUNTITLE{Data-Driven Inpatient Bed Assignment Using the $ P $ Model}

\TITLE{Data-Driven Inpatient Bed Assignment:
The $ P $ Model Approach to Boarding and Overflowing}

\ARTICLEAUTHORS{%
\AUTHOR{Shasha Han}
\AFF{Department of Analytics and Operations, NUS Business School, National University of Singapore, Singapore 117592\\ \EMAIL{shashahan@u.nus.edu}}
\AUTHOR{Shuangchi He}
\AFF{Department of Industrial Systems Engineering and Management, National University of Singapore, Singapore 117576\\ \EMAIL{heshuangchi@nus.edu.sg}}
\AUTHOR{Hong Choon Oh}
\AFF{Health Services Research, Eastern Health Alliance, Singapore 529541\\ \EMAIL{hong.choon.oh@easternhealth.sg}}
} 

\ABSTRACT{%

\emph{Problem definition}: Emergency department (ED) boarding refers to the practice of holding patients in the ED after they have been admitted to hospital wards, usually resulting from insufficient inpatient resources. Boarded patients may compete with new patients for medical resources in the ED, compromising the quality of emergency care. A common expedient for mitigating boarding is patient overflowing, i.e., sending patients to beds in other specialties or accommodation classes, which may compromise the quality of inpatient care and bring on operational challenges. We study inpatient bed assignment to shorten boarding times without excessive patient overflowing.

\emph{Academic/practical relevance}: As a cross-departmental issue, boarding is caused by periodic mismatches between the demand for inpatient care and the supply of inpatient resources. Besides reducing overall boarding times, mitigating the time-of-day effect is also essential for patients' safety. With numerous bed assignment constraints, it is challenging to strike a balance between boarding and overflowing.

\emph{Methodology}: We use a queue with multiple customer classes and multiple server pools to model hospital wards. Exploiting patient flow data from a hospital, we propose a computationally tractable approach to formulating the bed assignment problem, where the joint probability of all waiting patients meeting their respective delay targets is maximized.

\emph{Results}: By dynamically adjusting the overflow rate, the proposed approach is capable not only of reducing patients' waiting times, but also of mitigating the time-of-day effect on boarding times. In numerical experiments, our approach greatly outperforms both early discharge policies and threshold-based overflowing policies, which are commonly used in practice.

\emph{Managerial implications}: We provide a practicable approach to solving the bed assignment problem. This data-driven approach captures critical features of patient flow management, while the resulting optimization problem is practically solvable. The proposed approach is a useful tool for the control of queueing systems with time-sensitive service requirements.
}%


\KEYWORDS{emergency department boarding, healthcare operations, queueing network, routing, integer programming, $ P $ model}
\HISTORY{}

\maketitle

%


\section{Introduction}
\label{sec:introduction}

Emergency department (ED) crowding has been a worldwide crisis of healthcare delivery, compromising the quality of and access to emergency care in both developing and developed countries (\citealt{PinesETAL11b}). As a major cause of ED crowding, \emph{boarding} refers to the practice of holding patients in the ED after they have been admitted to hospital wards, usually resulting from insufficient inpatient resources. Because boarded patients may compete with new patients for medical resources in the ED, boarding could bring on treatment delays and higher complication rates. Boarded patients also suffer from prolonged waiting, since they may not receive the same level of care as in inpatient wards. Recent studies have identified an association between prolonged boarding times and increased mortality in both EDs and inpatient wards (\citealt{SingerETAL11,SunETAL13}). As a quality indicator of emergency care, boarding times are closely monitored by government agencies in some countries. In the United States, the Centers for Medicare and Medicaid Services include boarding times in the measures of their pay-for-performance program; on their Hospital Compare website, the statistics of boarding times from over 4,000 hospitals are accessible to the public. In Singapore, public hospitals must report boarding times to the Ministry of Health regularly; their daily median boarding times are published on the Ministry's website and updated every week. The Ministry may investigate those cases in which patients' waiting times for beds are excessively long.

\begin{figure}[t]
\centering
\begin{minipage}{.49\textwidth}
\includegraphics[trim={.3in 2.5in .4in 2.5in},height=2.35in]{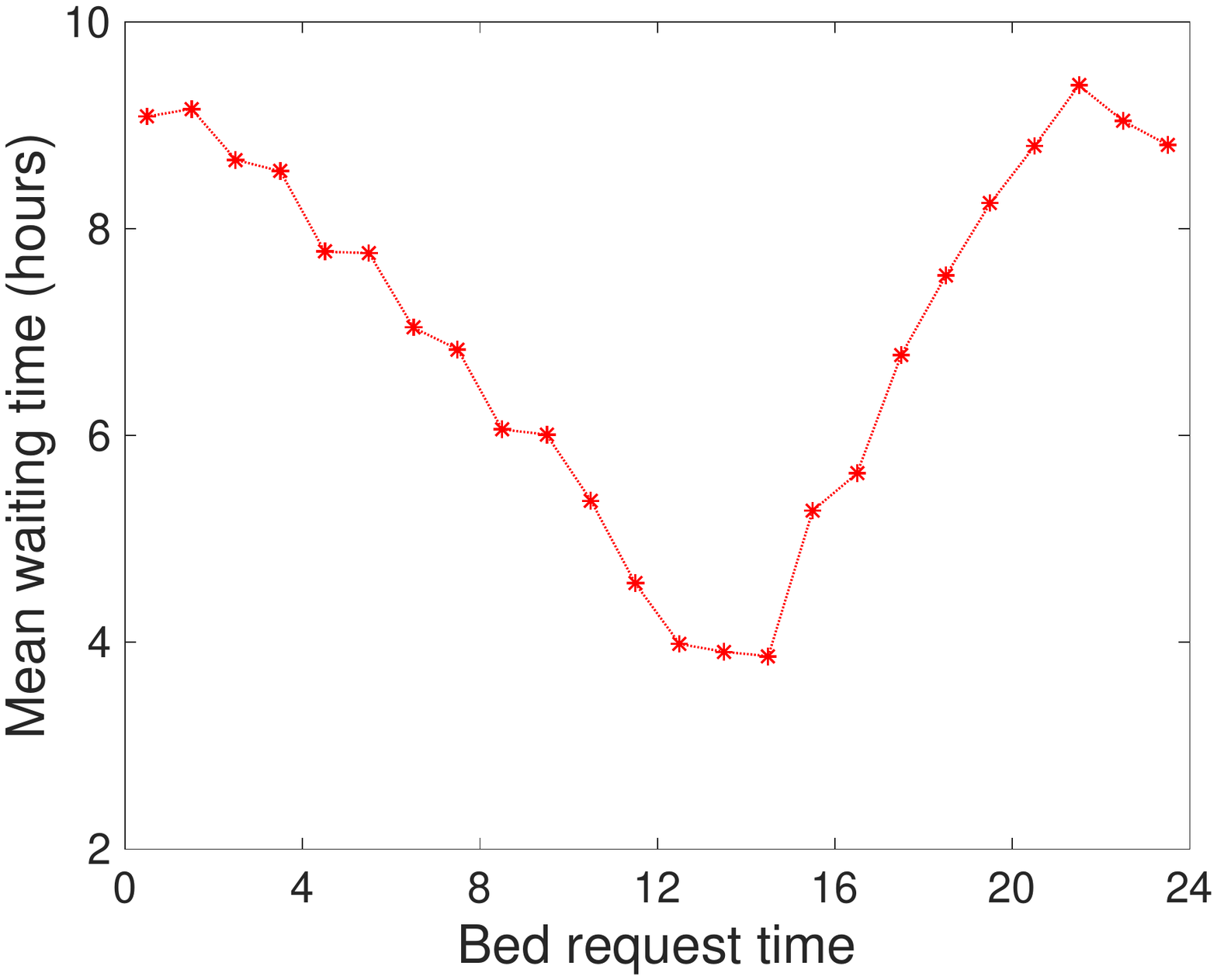}
\caption{Mean waiting times of ED patients at an anonymous public hospital.}
\label{fig:waiting-data}
\end{minipage}~~
\begin{minipage}{.49\textwidth}
\includegraphics[trim={.3in 2.5in .4in 2.5in},height=2.35in]{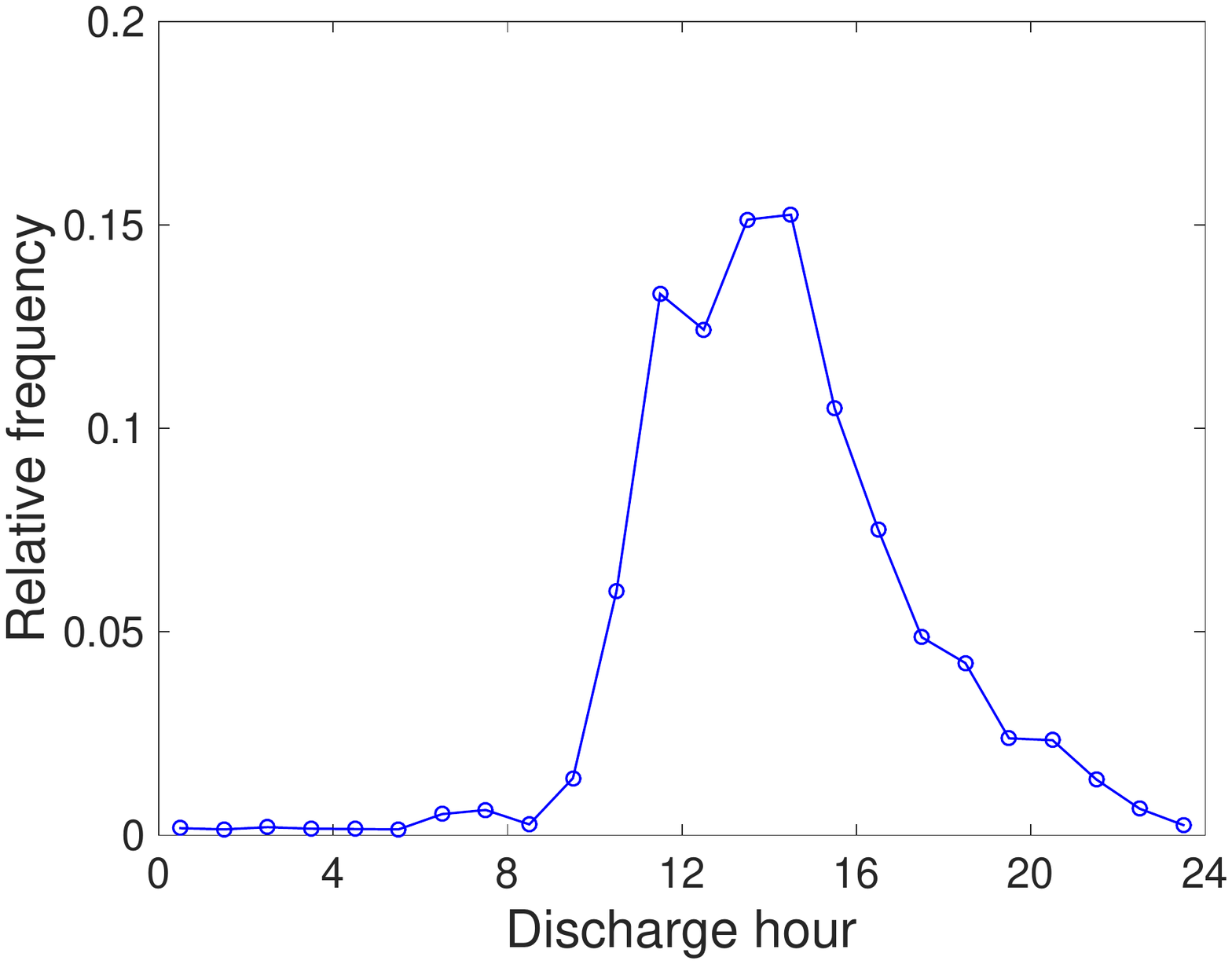}
\caption{Relative frequencies of inpatient discharge hours at the hospital.}
\label{fig:discharge-hours}
\end{minipage}
\end{figure}

As a cross-departmental issue, boarding is caused by \emph{periodic mismatches} between the demand for inpatient care and the supply of resources such as inpatient beds and hospital staff. In hospital wards, the majority of patients are commonly discharged in the afternoon, which may lead to a \emph{time-of-day effect} on boarding times: When the inpatient bed occupancy rate is high, patients who request beds at late night or in the early morning may experience long delays until appropriate beds are available. This phenomenon is evident in Figure~\ref{fig:waiting-data}, where the mean waiting times of ED patients at an anonymous public hospital in Singapore are plotted across different bed request hours. These results are estimated using the hospital's patient flow records of April--September 2015. Within this period, patients who requested beds at night were expected to wait more than eight hours, while those who requested beds in the early afternoon spent only about four hours. Prolonged boarding times could be more detrimental to ED patients who request beds at late night, because they may be in a serious condition which could be better managed in inpatient wards. Therefore, aside from reducing overall waiting times, mitigating the time-of-day effect on boarding times is also essential for patients' safety. We illustrate the relative frequencies of inpatient discharge hours at the hospital in Figure~\ref{fig:discharge-hours}. As one may expect, most patients were discharged in the afternoon, and not until then, would appropriate beds be available for a number of requests made on the previous night. \citet{ShiETAL16} reported a similar patient discharge pattern (see Figure~2 in their paper) and the resulting time-of-day effect on boarding times at another hospital in Singapore.

To mitigate ED boarding, some hospitals launch initiatives encouraging more patients to be discharged from wards before noon (\citealt{WertheimerETAL14,ShiETAL16}). Although positive outcomes are reported on patient throughput (\citealt{WertheimerETAL15}), there is no sufficient evidence in the literature that discharging patients earlier would be a practical solution to prolonged boarding times. To assess the influence of earlier discharge hours, \citet{ShiETAL16} developed a high-fidelity simulation model using a set of patient flow data from a public hospital in Singapore. They found out that in the hospital, boarding times may be equalized over the day if a portion of patients are discharged before 9\,\textsc{am}. However, such a discharge policy would be difficult to implement, because physicians and nurses are busy with ward rounds in the morning. Discharging patients at earlier hours will inevitably affect existing daily ward routines. As pointed out by \citet{Green08}, physicians are usually reluctant to adopt new operational practices owing to clinical concerns, and in general, hospital administrators are unwilling to exert pressure on physicians to change their behavior. In other words, resistance from medical staff may also be a hurdle in implementing a new patient discharge policy.

A common expedient for mitigating ED boarding is \emph{patient overflowing}, which refers to assigning patients to ``incorrect'' beds when the most appropriate beds are not available. For example, a cardiac patient could be sent to a bed in the general medical ward upon his admission, if there is no bed available in the cardiac ward. As a means of resource pooling, patient overflowing can always reduce excessive delays. However, it may compromise the quality of inpatient care, bringing on operational challenges such as additional costs and coordination between wards (\citealt{TeowETAL11}). Hence, a balance must be struck between boarding and overflowing.

We study inpatient bed assignment in this paper to shorten ED boarding times and mitigate the time-of-day effect. From a modeling perspective, hospital wards can be viewed as a queueing system, with patients being the customers and beds being the servers. Since both patients and beds are categorized by specialty, gender, and accommodation class, the assignment of inpatient beds turns out to be a routing problem for a queue with multiple customer classes and multiple server pools. We adopt the terminology used by \citet{ShiETAL16}, referring to a bed that completely matches a patient's specialty, gender, and requested accommodation class as a \emph{primary bed} for the patient, and referring to such a patient as a \emph{primary patient} for the bed. If a patient is allowed to be sent to a bed that does not match all the requirements, they are said to be a \emph{non-primary patient} and a \emph{non-primary bed}, respectively. With requirements for boarding times, it may not be possible to find every patient a primary bed. Sending a fraction of patients to non-primary beds would be necessary for meeting delay targets, especially when the bed occupancy rate is high. Since a non-primary bed for a boarded patient could be a primary bed for others, overflowing may prevent wards from accommodating more primary patients in the future. Without a carefully designed bed assignment policy, overflowing could be spread over multiple wards, resulting in more operational issues. \citet{ShiETAL16} reported that at a public hospital in Singapore, nearly 30\% of patients from the ED were sent to non-primary beds in 2008--2010.

To address this problem in a practical setting, we propose a dynamic bed assignment approach that exploits patient flow data from the hospital. When a bed request is received and there is no primary bed available, the assignment algorithm will determine whether to send the patient to a non-primary bed immediately, or to let the patient wait for a primary bed. When an inpatient bed becomes available and there are boarded patients, the algorithm will determine whether to assign a patient to this bed and who should be sent to this bed if there are multiple patients waiting. In other words, the assignment algorithm will determine \emph{whether to overflow} and \emph{where to overflow} for patients and beds.

We impose mandatory targets for patient boarding times on the bed assignment problem, in order to meet requirements for the quality of care. Our intention is to maximize the percentage of patients whose boarding times are within these targets, while the overflow rate is maintained below a certain level. We formulate an optimization problem to obtain assignment decisions by maximizing the joint probability of all boarded patients meeting the delay targets. This formulation allows us to shorten boarding times for bed requests at night, thus mitigating the time-of-day effect. Taking the joint probability of target attainment as the objective function, such an optimization problem was first studied by \citet{CharnesCooper63}, who termed the formulation the \emph{$ P $ model}. This formulation, however, has not been widely used in practice, in part because evaluating the joint probability requires integration in high dimensions, which is generally computationally difficult. In a recent study,  \citet{HeETAL18} pointed out that the $ P $ model could be tractable if the objective joint probability can be written into a product form. Considering a patient scheduling problem in EDs, they proposed a hybrid robust-stochastic framework under which the $ P $ model formulation can be modified into a computationally amiable form. They also demonstrated that their approach may outperform an asymptotically optimal scheduling approach proposed by \citet{HuangETAL15}.

Although both the present paper and the study by \citet{HeETAL18} concern control problems in queueing networks, the respective $ P $ model formulations differ greatly. First, the hybrid formulation by \citet{HeETAL18} imposes an additional constraint on the $ P $ model to convert the objective function into a product form, i.e., the feasible set of service times is required to be hyperrectangular. The optimal solution to their hybrid formulation is an approximate solution to the original problem. In the present paper, we obtain exact optimal solutions to the bed assignment problem, without using approximation techniques. This is because in the bed assignment problem, patients' waiting times are relatively short compared with their bed occupancy times. Such a queueing system operates in the \emph{Halfin--Whitt regime} (see, e.g., \citealt{HalfinWhitt81, BorstETAL04}), where the number of boarded patients is much smaller than the number of beds. Therefore, a patient's boarding time must be another patient's residual bed occupancy time. Since patients' discharge times within a day are assumed to be mutually independent, the joint probability of all boarded patients meeting the delay targets has a product form under any given assignment plan, which results in a tractable $ P $ model formulation. Second, the patient scheduling problem studied by \citet{HeETAL18} involves optimal sequencing, which is computationally demanding even if the objective joint probability has a product form. Only for small- to moderate-scale queueing networks, may the sequencing problem be practically solvable. In contrast, the $ P $ model can be used to solve routing problems for much larger queueing systems. In Section~\ref{sec:numerical}, we use the $ P $ model formulation to dynamically determine bed assignments for a hospital where patients are categorized into 50 types and the 571 inpatient beds are categorized into 34 pools. To the best of our knowledge, no other approaches in the literature are capable of solving routing problems in such a large queueing system with so many customer classes and server pools, when practical features such as general, yet non-identically distributed service times and time-sensitive service requirements are involved.

The main contribution of this paper is twofold. First, we provide a \emph{practicable} data-driven approach to solving the bed assignment problem in hospitals. Unlike conventional approaches in the literature relying on highly stylized models, the $ P $ model approach is able to capture critical features of patient flow management, such as a time-varying patient arrival process, general but non-identical bed occupancy times, and many patient and bed types, while the resulting optimization problem is still practically solvable. Using a set of patient flow data from a hospital, our simulation study shows that the proposed approach can greatly reduce patients' boarding times and mitigate the time-of-day effect. Second, through this bed assignment problem, we demonstrate that the $ P $ model could be a useful tool for the control of queueing systems with time-sensitive service requirements. Solving such problems is usually difficult under the conventional probabilistic framework of queueing theory. We expect to see more applications of the $ P $ model in large-scale stochastic systems such as call centers and computer clusters.

The remainder of this paper is organized as follows. By reviewing related literature, we position our study in Section~\ref{sec:literature}. We introduce the multi-class, multi-pool queueing model for the bed assignment problem in Section~\ref{sec:model}. In Section~\ref{sec:P-model}, we present a tractable approach to solving this problem, using the $ P $ model formulation. We introduce a high-fidelity simulation model in Section~\ref{sec:simulation}, which is primarily based on the simulation model by \citet{ShiETAL16}. Using the simulation model, we conduct a data-based numerical study in Section~\ref{sec:numerical} to assess the $ P $ model approach. The paper is concluded in Section~\ref{sec:conclusion}. We leave the proof of a proposition and details about the simulation model in the appendix.

\section{Related Literature}
\label{sec:literature}

The assignment of inpatient beds is a routing problem for a queueing system with many servers, where most customers should be sent to matched servers. In the literature, routing problems for many-server queues are primarily motived by call center applications. Optimal call routing policies have been studied by \citet{BassambooETAL06, DaiTezcan08, BassambooZeevi09, GurvichWhitt09a, GurvichWhitt09b, GurvichWhitt10, ArmonyMandelbaum11}, and \citet{StolyarTezcan11} under various operational regimes and cost structures. Most of these studies rely on asymptotic analysis and Markovian assumptions such as independent, exponentially distributed service times. The bed assignment problem is generally more difficult, because some assumptions in call routing are no longer applicable in the hospital setting. For example, \citet{ShiETAL16} pointed out that the discharge pattern in hospital wards may affect boarding times significantly. As a consequence, one should not assume identically distributed bed occupancy times (see \citealt{DaiShi17,CDG17} for more discussion). With constraints on specialty, gender, and accommodation class, there may be tens of patient types and bed pools, while some pools may have only several beds. Even though there are hundreds of beds in a hospital, the routing policies that are asymptotically optimal for large queues may not work well for all patient types and bed pools. Therefore, we cannot adopt routing policies for call centers to solve the bed assignment problem. One may refer to \citet{ArmonyETAL15} for an overview of the general patient transfer process to hospital wards.

As a major concern of inpatient bed assignment, how to balance excessive delays against overflowing has attracted considerable attention from operations researchers. Using a simulation model calibrated with patient flow data, \citet{HarrisonETAL05} analyzed bed occupancy rates in a hospital in Australia, pointing out that both the seasonal variations of patient arrival rates and the variable patient discharge rates over days contributed a lot to frequent overflowing in the hospital. \citet{ThompsonETAL09} studied periodic reallocation of patients across multiple hospital wards to reduce waiting times for beds. They modeled this problem as a finite-horizon Markov decision process (MDP), using random sampling to simplify computation. It was reported that this patient reallocation approach helped a hospital in Connecticut to cut the average waiting time for beds by half and to increase the hospital's revenue by 1\%. Transferring patients between wards during their inpatient stays, however, may add to staff's workload and interfere with the wards' daily routines. It may not be a viable option for all hospitals. \citet{SamiedaluieETAL17} studied a patient admission problem within a neurology ward. They also used an MDP to determine whether an incoming patient should be admitted to the ward and which patient should be sent to each available bed. They assumed that patients would be transferred to other hospitals if not admitted to the ward. Patient overflowing across wards is not an option in their paper.

Using a high-fidelity simulation model, \citet{ShiETAL16} examined the bed assignment policy employed by a public hospital in Singapore, where overflowing was triggered when a patient's waiting time exceeded a predetermined threshold. They also analyzed the factors that may cause prolonged delays and proposed early discharge policies to mitigate the time-of-day effect on boarding times. By analyzing a simplified queueing model with two customer classes and two server pools, \citet{KilincETAL16} found some structural properties of optimal bed assignment policies. They proved that under such a policy, a patient will be sent to a non-primary bed only when his waiting time exceeds a threshold that depends on the system's current state. Based on their findings, the authors proposed heuristic policies that are implementable in hospitals.

An independent, concurrent study by \citet{DaiShi18} is the most relevant work to our paper. They adopted a two-time-scale model proposed by \citet{ShiETAL16} for bed occupancy times, and formulated the bed assignment problem as a discrete-time, infinite-horizon MDP. Their formulation allows for practical features such as time-varying patient arrival and discharge patterns, thus having a high-dimensional state space. To cope with the curse of dimensionality, they resorted to the technique of approximate dynamic programming (ADP). Although our study also follows the two-time-scale model to formulate bed occupancy times, the approach proposed in this paper is methodologically distinct, not relying on the MDP framework. By solving the $ P $ model problem repeatedly, we provide a computationally efficient approach that is readily implementable for hospitals with various constraints on their inpatient management. 

We would point out that our dynamic $ P $ model approach does not satisfy Bellman's principle of optimality. This approach is \emph{myopic} in nature, although future bed request rates could be incorporated in the formulation (through the overflow budget given by \eqref{eq:overflow-budget}). We would like to explore such an approach because patient arrival patterns in hospitals are highly variable. In an uncertain environment, a computationally tractable myopic policy may perform well, because its simple structure allows us to incorporate some realistic features that are difficult to be included in the MDP framework (see \citealt{HeETAL18} for more discussion). In Section~\ref{sec:P-model}, we translate the $ P $ model formulation into a tractable integer linear program, the solution to which can be quickly obtained using standard solvers.

\section{A Multi-Class, Multi-Pool Queueing Model}
\label{sec:model}

We use a multi-class, multi-pool queueing system to model hospital wards, which are managed by a centralized bed assignment system to meet requirements for boarding times without sending excessive patients to non-primary beds. Based on the state of boarded patients and the availability of beds, the bed assignment system will make sequential recommendations for incoming patients and available beds.

In our partner hospital, there is a department in charge of inpatient bed management, called the Bed Management Unit (BMU). The general bed assignment process in the hospital is as follows. When a bed request is received by the BMU, a staff member will first search the wards for a primary bed that matches the patient's specialty, gender, and requested accommodation class. If a primary bed is available, the BMU will assign the patient to this bed immediately and initiate the patient transfer process. Otherwise, the patient may have to wait until a primary bed is available. If a patient has been waiting for a long time or no primary bed is expected to be available within hours, the BMU may assign the patient to an available non-primary bed.

\begin{figure}[t]
\centering
\includegraphics[height=2.2in]{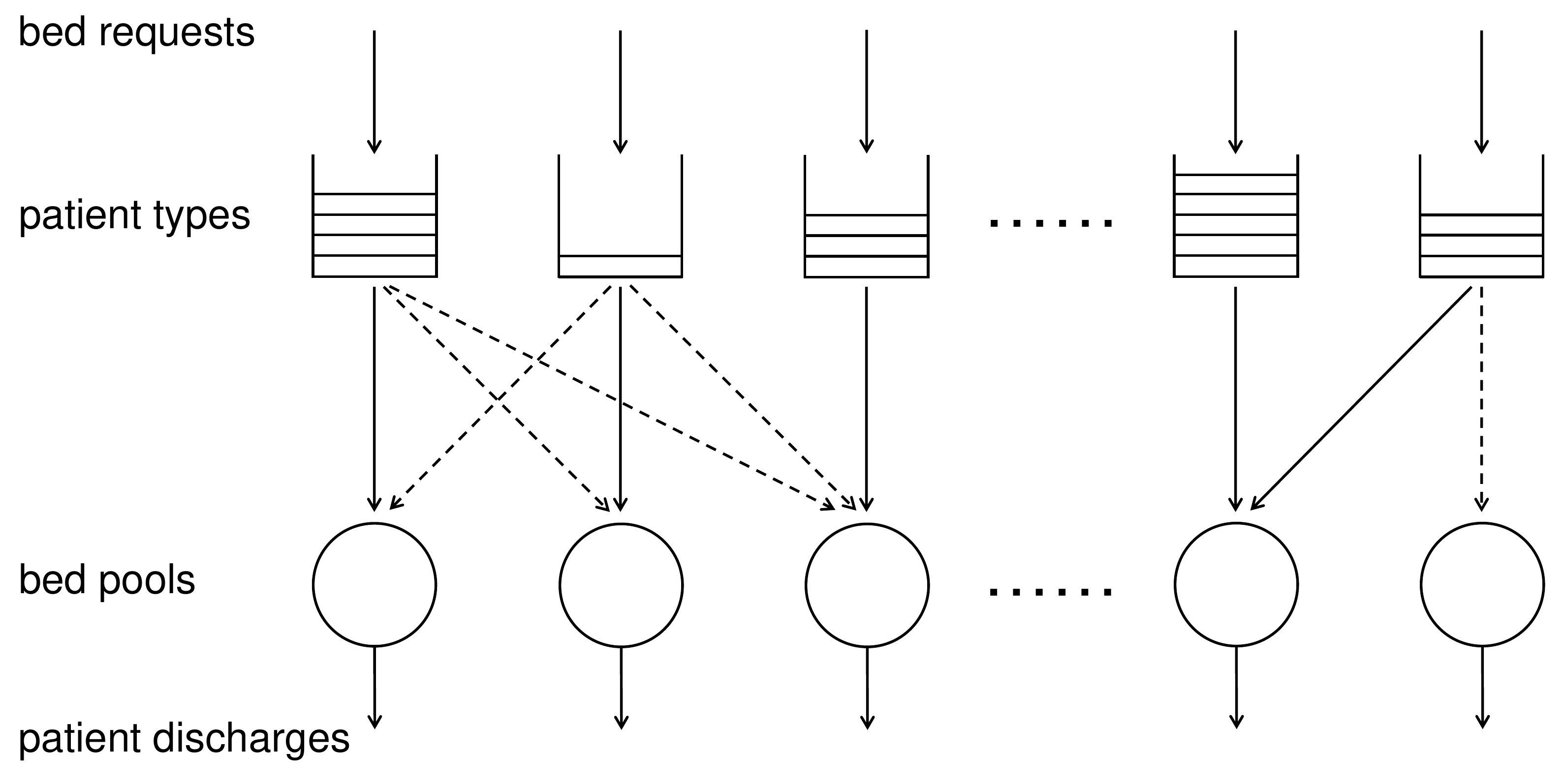}
\medskip
\caption{The multi-class, multi-pool queueing model for hospital wards. Solid lines represent primary patient-bed pairs and dashed lines represent non-primary pairs.}
\label{fig:queueing-model}
\end{figure}

Both patients and beds are categorized by specialty, gender, and accommodation class. In other words, this queueing system has multiple customer types and multiple server pools. The gender constraint requires that each inpatient room must accommodate patients of the same gender, i.e., while a single room is gender-neutral, a room with multiple beds can only accept patients of the same gender as the current patients. As a result, each bed is labeled as \emph{male only}, \emph{female only}, or \emph{gender-neutral}. With this requirement, some beds may have more than one type of primary patients, e.g., a single room in the surgical ward (which has a gender-neutral bed) may accept either a male or female surgical patient. Therefore, the number of patient types could be different from the number of bed pools in this queueing model. We allow a bed pool to have more than one primary patient type, and a patient type to have more than one primary bed pool. This assumption makes our queueing model slightly different from the models in \citet{KilincETAL16} and \citet{DaiShi18}, where each bed pool is assumed to have one and only one primary patient type.

We assume that bed requests are sent to the BMU in a stochastic manner. The bed request process may have a time-varying rate. The BMU must process all requests. In other words, the queueing system has infinite buffers holding all waiting customers. We illustrate this multi-class, multi-pool queueing model in Figure~\ref{fig:queueing-model}. We assume that each patient is required to be transferred to an inpatient bed within a given delay target. To meet this requirement, patient overflowing is generally inevitable. However, the BMU must also ensure that most patients are sent to primary beds so that the overall quality of care will not be compromised. An optimization model should thus be formulated in order for the hospital to strike a balance between reducing excessive delays and matching patients with primary beds.

\section{The $ \boldsymbol{P} $ Model Approach}
\label{sec:P-model}

We will design a dynamic bed assignment algorithm driven by patient flow data. Our intention is to maximize the percentage of patients whose boarding times do not exceed their respective delay targets, while maintaining the overflow rate at a reasonable level. We assume that a decision iteration is triggered at the following event times: (i) when a bed becomes available, (ii) when a bed request is received and there is at least one available bed for the patient, (iii) when a patient's boarding time is about to exceed the delay target. Each time, the bed assignment algorithm will recommend the next patient to be assigned to an available bed. 

Let $ t $ be a moment when the bed assignment algorithm is triggered. We use $ \mathcal{I} $ to denote the set of patients whose bed requests are pending at time $ t $ and $ \mathcal{J} $ the set of beds that are either available at time $ t $ or to be available soon (see more discussion later). For the sake of convenience, the dependence on $ t $ is suppressed in notation. We assume that both $ \mathcal{I} $ and $ \mathcal{J} $ are nonempty. For $ i \in \mathcal{I} $, we use $ \mathcal{J}_{i} \subset \mathcal{J} $ to denote the set of primary and non-primary beds for patient~$ i $. Without loss of generality, we assume that $ \mathcal{J} = \bigcup_{i \in \mathcal{I}} \mathcal{J}_{i}$.

For $ j \in \mathcal{J} $, let $ \tilde{d}_{j} $ be the time when bed~$ j $ becomes available. If bed~$ j $ is available at time $ t $, we simply take $ \tilde{d}_{j} = t $. If $ \tilde{d}_{j} > t $, i.e., bed~$ j $ will be available later, $ \tilde{d}_{j} $ is a random variable whose distribution can be estimated from data, e.g., using the discharge time distribution illustrated in Figure~\ref{fig:discharge-hours}. We assume that $ \{ \tilde{d}_{j}:j\in\mathcal{J} \} $ is a set of independent random variables. (This assumption can be verified by a standard $ \chi^{2} $ test of independence. It is generally valid for hospitals having many inpatient wards.) Let $ a_{i} $ be patient $ i $'s bed request time, with $ a_{i} \leq t $ for all $ i \in \mathcal{I} $. Patient~$ i $ is allowed to be assigned to bed~$ j $ only if $ j \in \mathcal{J}_{i} $, in which case the boarding time of patient~$ i $ will be $ \tilde{d}_{j} - a_{i} $. We use $ \tau_{i} $ to denote the delay target of patient~$ i $, assuming $ \tau_{i} > 0 $ for all $ i \in \mathcal{I} $. We would find a bed assignment plan for the boarded patients to maximize the joint probability of all of them meeting their respective delay targets.

We use a binary variable $ z_{ij} \in \{0,1\} $ to indicate the assignment decision of patient~$ i $ to bed~$ j $, i.e., $ z_{ij} = 1 $ if he is assigned to the bed and $ z_{ij} = 0 $ otherwise. An assignment plan $ \boldsymbol{z} = (z_{ij})_{i \in \mathcal{I},j\in \mathcal{J}} $ is said to be \emph{admissible} if every waiting patient is assigned to exactly one bed that is primary or non-primary for him, i.e., 
\[
\sum_{j \in \mathcal{J}_{i}} z_{ij}  = 1 \quad \mbox{for } i \in \mathcal{I},
\]
and each bed is allowed to accommodate at most one patient, i.e.,
\[
\sum_{i\in\mathcal{I}} z_{ij} \leq 1 \quad \mbox{for } j \in \mathcal{J}.
\]
We use $ \mathcal{A} $ to denote the set of all admissible assignment plans. We assume that $ \mathcal{J} $ includes sufficient beds so that $ \mathcal{A} $ is nonempty. Typically, we may take $ \mathcal{J} $ to be the set of beds that are either available now or to be available on the present day. If $ \mathcal{A} $ is empty with such a choice of $ \mathcal{J} $, we may include more beds that will be available later (e.g., all beds to be available the next day) to make $ \mathcal{A} $ nonempty. This assumption is reasonable because in general, a patient's discharge is scheduled at least one day in advance. Under an arbitrary assignment plan $ \boldsymbol{z} \in \mathcal{A} $, the joint probability of all waiting patients meeting their delay targets is given by
\[  
\mathbb{P}\big( \tilde{d}_{j} - a_{i} \leq \tau_{i} : z_{ij} = 1,\ i \in \mathcal{I},\ j \in \mathcal{J}_{i}\big) = \prod_{i\in\mathcal{I}} \prod_{j\in\mathcal{J}_{i}} \mathbb{P}\big(\tilde{d}_{j} - a_{i} \leq \tau_{i}\big) ^ {z_{ij}},
\]
where the product form on the right side follows from the independence assumption of $ \{ \tilde{d}_{j}:j\in\mathcal{J} \} $.

We need to send most patients to primary beds. To this end, a cost is incurred when a patient is assigned to a non-primary bed. Let $ u_{ij} \geq 0 $ be the cost of assigning patient $ i \in \mathcal{I}$ to bed $ j \in \mathcal{J}_{i} $. We assume that $ u_{ij} = 0 $ if bed~$ j $ is a primary bed for patient~$ i $, and that $ u_{ij} > 0 $ if bed~$ j $ is a non-primary bed. To refrain from sending too many patients to non-primary beds, we add an extra constraint, requiring the total overflow cost not to exceed an \emph{overflow budget} $ B > 0$. Hence, we may solve the following $ P $ model problem to find an optimal assignment plan:
\begin{equation}
\begin{array}{l@{\quad}l@{}l@{\quad}l@{}}
\text{max} & \multicolumn{3}{l}{\displaystyle{\sum_{i\in\mathcal{I}} \sum_{j\in\mathcal{J}_{i}} z_{ij}\ln \mathbb{P}\big(\tilde{d}_{j} - a_{i} \leq \tau_{i}\big)}}\\
\text{s.t.} & \displaystyle{\sum_{j \in \mathcal{J}_{i}} z_{ij}  = 1}, &  & i \in \mathcal{I} \\
& \displaystyle{\sum_{i\in\mathcal{I}} z_{ij} \leq 1}, &  & j \in \mathcal{J}\\
& \displaystyle{\sum_{i\in\mathcal{I}} \sum_{j\in\mathcal{J}_{i}} z_{ij}u_{ij} \leq B}, & &\\
& z_{ij} \in \{0,1\},& & i \in \mathcal{I},\ j \in \mathcal{J}.
\end{array}
\label{eq:optimal-assignment}
\end{equation}
In this formulation, patient overflowing is controlled by the overflow budget $ B $. We may take $ u_{ij} = 1 $ if bed~$ j $ is a non-primary bed for patient~$ i $. In this case, the overflow budget turns out to be the maximum number of non-primary assignments allowed by the current iteration. By adjusting the overflow budget, the overflow rate of patients can be maintained at a reasonable level.

A greater overflow budget allows more patients to be assigned to non-primary beds, thus yielding shorter boarding times. However, a greater overflow budget may also result in worse quality of care, so $ B $ should be kept as small as possible. In order to get a feasible solution to \eqref{eq:optimal-assignment}, the minimum overflow budget required could be obtained by solving the following integer program:
\begin{equation}
\begin{array}{ll@{\quad}l@{}l@{\quad}l@{}}
B^{\star} = & \text{min} & \displaystyle{\sum_{i\in\mathcal{I}} \sum_{j \in \mathcal{J}_{i} } x_{ij} u_{ij} }\\
& \text{s.t.} & \displaystyle{\sum_{j\in\mathcal{J}_{i}} x_{ij} = 1}, & & i\in \mathcal{I}\\
& & \displaystyle{\sum_{i\in\mathcal{I}} x_{ij} \leq 1}, & & j\in\mathcal{J}\\
& & x_{ij} \in \{0,1\}, & & i \in \mathcal{I},\ j \in \mathcal{J},
\end{array}
\label{eq:feasibility}
\end{equation}
which is an \emph{assignment problem}. This problem is efficiently solvable because the linear programming (LP) relaxation of that has integral solutions. We may thus solve the LP relaxation instead (see, e.g., \citealp{Balinski86}). More specifically, in the multi-class, multi-pool queueing model, all compatible customer-server pairs form a bipartite graph (see Figure~\ref{fig:queueing-model}). Since the node-edge incidence matrix of a bipartite must be totally unimodular, routing constraints in the queueing model specify an integral polyhedron in the LP relaxation (see, e.g., Corollary~III.1.2.9 in \citealp{NemhauserWolsey88}). These constraints also appear in the integer linear program \eqref{eq:optimal-assignment}, the feasible region of which is nonempty if $ B \geq B^{\star} $. Although the feasible region of \eqref{eq:optimal-assignment} may not be an integral polyhedron with the extra constraint on the overflow cost, this formulation may still have a ``nice'' structure, allowing us to quickly solve \eqref{eq:optimal-assignment} using existing integer programming solvers. A computational comparison is made between the $ P $ model formulation \eqref{eq:optimal-assignment} and the LP relaxation in Section~\ref{sec:computational}. We find out that it is generally efficient to obtain an optimal solution to the $ P $ model formulation.

When a bed becomes available, the bed assignment algorithm determines the next patient for this bed by solving \eqref{eq:optimal-assignment}. The following proposition states that there exists an optimal assignment plan under which this bed will take a patient whose bed request time is the earliest in his type. For $ i,k \in \mathcal{I} $, we say patients~$ i $ and~$ k $ are of \emph{the same type} if $ \mathcal{J}_{i} = \mathcal{J}_{k}$, $ \tau_{i} = \tau_{k} $, and $ u_{ij} = u_{kj} $ for all $ j \in \mathcal{J}_{i} $. The proof of Proposition~\ref{prop:FCFS} can be found in Section~\ref{sec:proofs} in the online appendix.

\begin{proposition}
\label{prop:FCFS}
Assume that $ B \geq B^{\star} $ in \eqref{eq:optimal-assignment}, where $ B^{\star} $ is the minimum overflow budget given by \eqref{eq:feasibility}. For an arbitrary bed $ j \in \mathcal{J} $ that is available to take the next patient, there exists an optimal solution $ \boldsymbol{z}^{\star} $ to \eqref{eq:optimal-assignment} that satisfies the following condition: If patient~$ k $ is assigned to bed~$ j $ under $ \boldsymbol{z}^{\star} $, then his bed request time must be the earliest in his type, i.e., $ a_{k} \leq a_{\ell} $ for any patient~$ \ell $ of the same type as patient~$ k $.
\end{proposition}

Assume that bed~$ j $ is available and patient~$ k $ is assigned to bed~$ j $ according to an optimal solution to \eqref{eq:optimal-assignment}. If patient~$ k $ does not have the earliest bed request time in his type, we may find the earliest one and send that patient to bed~$ j $ instead. As shown in the proof of Proposition~\ref{prop:FCFS}, switching beds between these two patients will not affect the optimality of the assignment plan. In this way, we can obtain an optimal assignment plan that satisfies the condition in Proposition~\ref{prop:FCFS} based on an arbitrary optimal solution to \eqref{eq:optimal-assignment}.

We determine bed assignments for incoming patients by solving \eqref{eq:optimal-assignment} repeatedly. In this process, a decision iteration is triggered when a bed becomes available, when a bed request arrives with at least one bed being available, or when a patient's boarding time is about to exceed the delay target. In general, patients are not sent to beds on a first-come, first-served (FCFS) basis, because different patients may have different urgency levels and delay requirements. However, if we follow an optimal assignment plan that satisfies the condition in Proposition~\ref{prop:FCFS} in each iteration, patients belonging to each type will be sent to beds in the FCFS manner. This feature of the $ P $ model formulation not only guarantees fairness in serving patients of the same type, but also prevents excessively long waiting times patients might experience when their delay targets cannot be met.

Even if there are available beds when a decision iteration is triggered, it is still possible that by an optimal solution to \eqref{eq:optimal-assignment}, all patients should be assigned to beds to be available later. This may happen when all patients' boarding times are far below their delay targets, so it is not necessary to send anyone to a non-primary bed immediately. Such a phenomenon implies that our dynamic bed assignment policy is \emph{not} non-idling. As pointed out by \citet{KilincETAL16} and \citet{DaiShi18}, it is necessary to idle some beds from time to time so that the overflow rate can be maintained at a reasonable level. The practice of strategic idling is also helpful for improving operational performance in other service systems, e.g., \citet{BaronETAL17} studied idling policies in a health examination network.

A decision iteration is triggered when a bed request is received and there is at least one bed available for the patient. If a primary bed is available, one may expect the incoming patient to be sent to this bed immediately. However, the $ P $ model approach may not always make such a decision. An optimal solution to \eqref{eq:optimal-assignment} occasionally specifies a non-primary patient for the available bed, even if the incoming primary patient will be waiting. This may happen when the non-primary patient's boarding time is about to exceed the delay target, in which case it is necessary to assign the patient to a non-primary bed immediately. Although such a decision is important in preventing a patient's boarding time from exceeding the delay target, it will not happen under the ADP-based policy by \citet{DaiShi18}, who assumed that if a primary bed is available upon a bed request, the incoming patient will always be assigned to the primary bed immediately. Their assumption is reasonable because no delay targets are imposed on their formulation.

To ensure that an optimal solution to \eqref{eq:optimal-assignment} exists, we should take $ B \geq B^{\star}$ in each decision iteration. The selection of the overflow budget will affect the solution to \eqref{eq:optimal-assignment}, thus having significant influence on the performance of the bed assignment algorithm. We adjust this parameter dynamically to mitigate the time-of-day effect. More specifically, if a decision iteration is triggered at time $ t $, we determine the overflow budget for this iteration by the following formula:
\begin{equation}
\label{eq:overflow-budget}
B = \max \Big\{ B^{\star},\ \alpha \cdot |\mathcal{I}| + \beta \int_{t}^{t+\Delta}\lambda (u)\,\mathrm{d}u \Big\},
\end{equation} 
where $ B^{\star} $ is the minimum overflow budget given by \eqref{eq:feasibility}, $ \lambda(u) $ is the patient arrival rate at time $ u $, and $ \alpha $, $ \beta $, and $ \Delta $ are three positive numbers to be determined heuristically. Clearly, the overflow budget given by \eqref{eq:overflow-budget} satisfies $ B \geq B^{\star} $. It takes such a form because of the following considerations. First, in order to mitigate the time-of-day effect, we should impose a greater overflow budget when the system is more congested. The term $ \alpha \cdot |\mathcal{I}| $ in \eqref{eq:overflow-budget} allows us to adjust $ B $ based on the number of waiting patients. Second, if bed requests are expected to surge in the next few hours, sending more patients to non-primary beds in advance would be helpful in preventing congestion. Hence, we also adjust the overflow budget using an estimate of the number of bed requests in the next $ \Delta $ hours, which is the integral term in \eqref{eq:overflow-budget}. When all $ u_{ij} $'s are integer-valued, so is $ B^{\star} $. The problem \eqref{eq:optimal-assignment} can be solved in a more efficient way by integer programming solvers if $ B $ is also an integer. In this case, we would modify \eqref{eq:overflow-budget} into 
  \begin{equation}
\label{eq:overflow-budget-integer}
B = \max \Big\{ B^{\star},\ \Big\lfloor \alpha \cdot |\mathcal{I}| + \beta \int_{t}^{t+\Delta}\lambda (u)\,\mathrm{d}u \Big\rfloor \Big\}
\end{equation} 
for the convenience of computation.

\section{Simulation Model}
\label{sec:simulation}

In this section, we introduce the simulation model for inpatient operations in our partner hospital. This simulation model is primarily based on the stochastic network model proposed by \citet{ShiETAL16}, but including several features absent from their model, such as gender and accommodation class constraints. This simulation platform is populated with parameters estimated from patient flow records, allowing us to faithfully reproduce patient flow dynamics in this hospital. We replicate important performance measures from the data set to validate this model, including ED patients' mean waiting times, service levels, and overflow rates. This simulation model will be used for assessing the proposed bed assignment approach in Section~\ref{sec:numerical}.

\subsection{Categorization of Patients and Beds}
\label{sec:environment}

\begin{figure}[t]
\centering
\begin{minipage}{0.49\textwidth}
\includegraphics[trim={.3in 2.5in .4in 2.5in},height=2.35in]{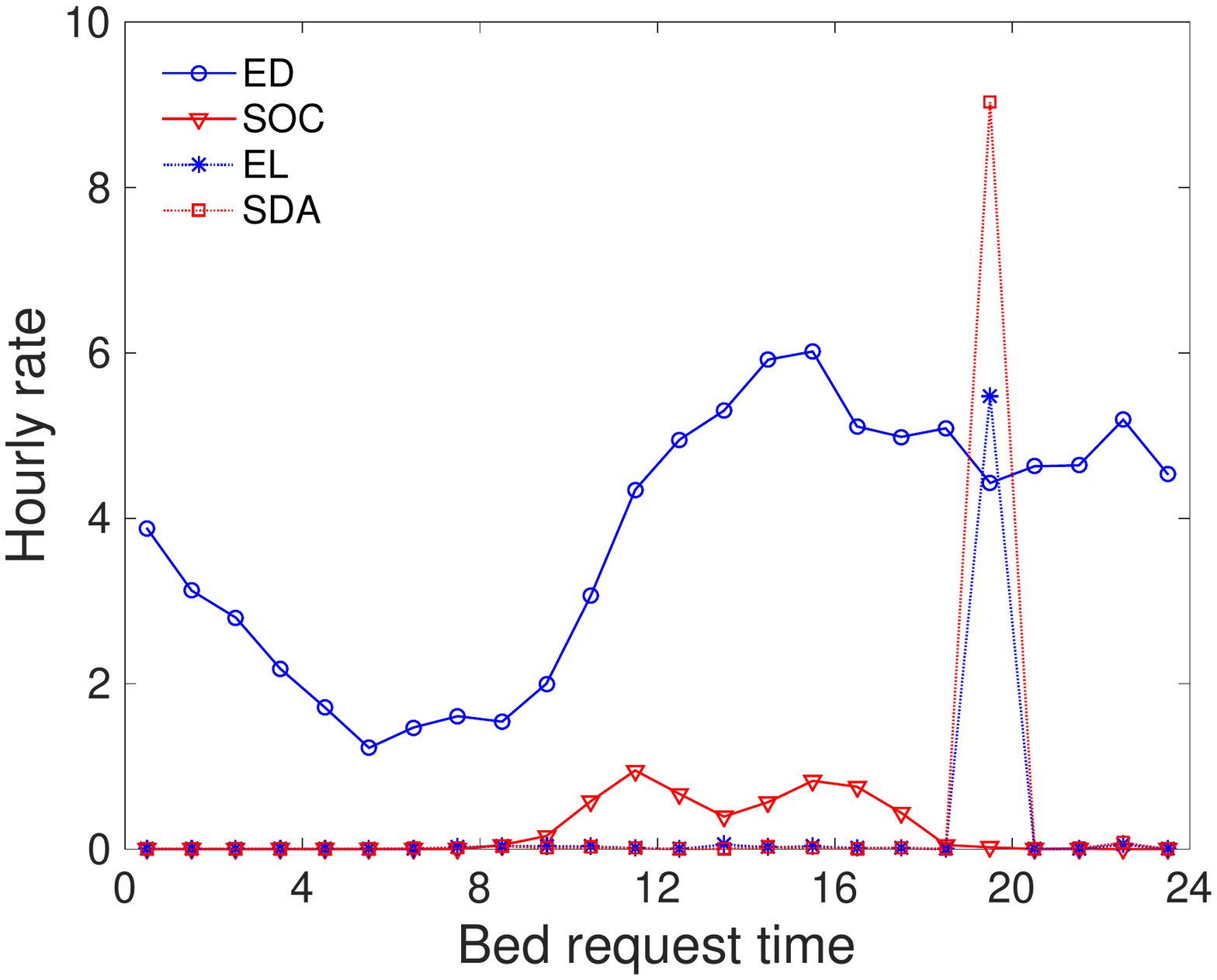}
\caption{Hourly bed request rates from four sources.}
\label{fig:arrival-rates}
\end{minipage}~~
\begin{minipage}{.49\textwidth}
\includegraphics[trim={.3in 2.5in .4in 2.5in},height=2.35in]{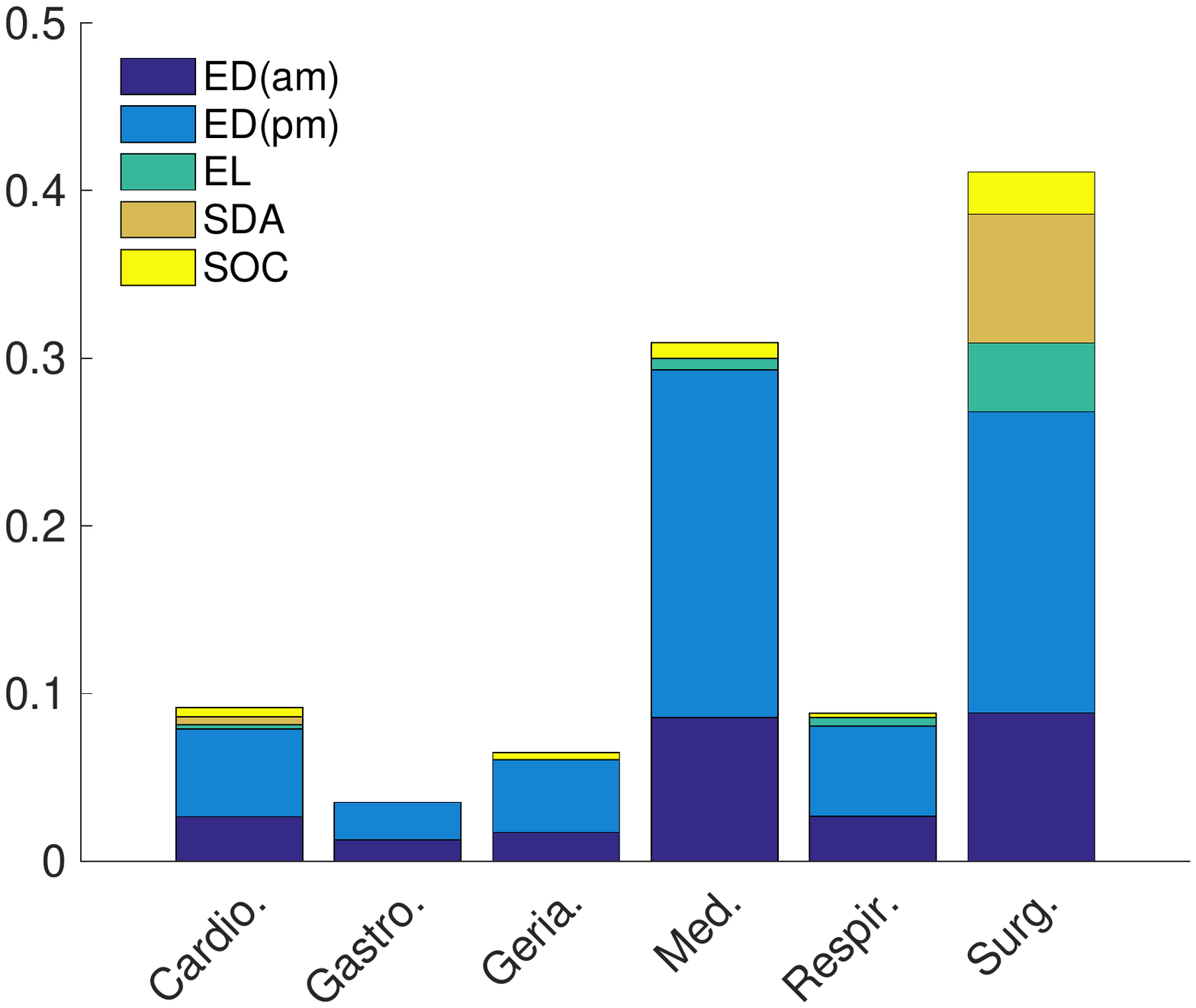}
\caption{Bed request proportions of six specialties.}
\label{fig:specialties}
\end{minipage}
\end{figure}

The hospital provided us with their patient flow records of April--September 2015. On average, the hospital wards admitted about 134 patients per day from four different sources: 81.4\% of inpatient admissions were requested from the ED, 8.4\% were same-day-admissions (SDA) requested for patients having day surgery, 5.2\% were elective (EL) admissions scheduled in advance, and 5.0\% were requested from specialist outpatient clinics (SOC). The hourly bed request rates from the four sources are plotted in Figure~\ref{fig:arrival-rates}, where those from ED and SOC patients are time-varying throughout the day and almost all bed requests from EL and SDA patients were sent to the BMU around 7\,\textsc{pm}. According to the hospital's weekly schedule, there are no bed requests from EL patients on Friday and no bed requests from SDA patients on Saturday. In the simulation model, we assume that bed requests from ED and SOC patients follow two periodic, non-homogeneous Poisson processes with the period being one day, while bed requests from EL and SDA patients are received at 7\,\textsc{pm}, except on Friday and Saturday, respectively.

All inpatient beds in the hospital are broadly categorized into six specialties, i.e., cardiology, gastroenterology, geriatrics, general medicine, respiratory medicine, and surgery. The proportions of bed requests from these specialties are illustrated in Figure~\ref{fig:specialties}, where the fractions from different sources are also specified. We differentiate ED patients who were admitted before and after 12\,\textsc{pm}, because these two groups of patients have different length-of-stay distributions. This phenomenon was discussed by \citet{ShiETAL16}, who speculated that ED patients admitted in the morning could undergo some diagnostic tests in the afternoon, whereas those admitted later may have to wait until the next day. A fraction of patients would thus spend one more day in the hospital if admitted after 12\,\textsc{pm}. We report the mean lengths of stay and the corresponding standard deviations of the six specialties from different patient sources in Table~\ref{tab:LOS} in Section~\ref{sec:more-simulation}. The mean lengths of stay of ED patients admitted in the morning are shorter in all specialties.

In this hospital, inpatient rooms are categorized into four accommodation classes, A1, B1, B2, and C, having one, four, six, and eight beds, respectively. The percentages of bed requests in the six specialties for these four classes are reported in Table~\ref{tab:accommodation} in Section~\ref{sec:more-simulation}. The beds of classes B2 and C are subsidized, and thus demand for these two classes is much higher than the other two classes. As we mentioned earlier in Section~\ref{sec:model}, each room must accommodate patients of the same gender. While A1 beds are gender-neutral, all other beds can only accommodate patients of a particular gender. The gender distributions in the six specialties are reported in Table~\ref{tab:gender} in Section~\ref{sec:more-simulation}. We would emphasize that patient overflowing refers to sending patients not only to wards in other specialties, but also to beds belonging to the requested specialties but in other accommodation classes.

According to the patient flow records, the number of beds in several wards varied across the period of interest owing to renovation or maintenance, opening of new wards, etc. For the convenience of implementation, we assume that the number of beds is fixed in each ward in simulation. To match the performance statistics of the data, we have slightly adjusted the number of beds in a few wards. There are totally 571 inpatient beds in the simulation model, divided into 34 pools by specialty, gender, and accommodation class. Please refer to Table~\ref{tab:bed-pools} in Section~\ref{sec:more-simulation} for the indexes and capacities of these bed pools. All patients are categorized into 50 types accordingly. The primary and non-primary bed pools for male and female patients are specified in Tables~\ref{tab:overflow-male} and~\ref{tab:overflow-female} in Section~\ref{sec:more-simulation}, respectively, according to the hospital's guidelines. In these two tables, both \emph{preferred} and \emph{secondary} pools refer to non-primary bed pools for a particular type of patients. When a patient has to be sent to a non-primary bed, the BMU will try their best to find the patient a preferred bed because the conditions are more similar to a primary bed. Only if a preferred bed is not available, would the BMU assign the patient to a secondary bed. The overflow rate refers to the percentage of patients sent to either preferred or secondary beds.

\subsection{Bed Occupancy Times and Pre- and Post-Allocation Delays}
\label{sec:two-time-scale}

\citet{ShiETAL16} pointed out that a patient's bed occupancy time generally follows a two-time-scale model, i.e.,
\begin{equation}
\label{eq:two-time-scale}
T = L - A + D,
\end{equation}
where $ L $ is the length of stay or number of nights the patient stays at the ward, $ A $ is the time on the admission day when the patient is admitted to the bed, and $ D $ is the time on the discharge day when the bed is vacated. Therefore, if measured in days, $ L $ must be an integer, $ 0 \leq A < 1 $, and $ 0 \leq D < 1 $. Although it is reasonable to assume that both lengths of stay and discharge times are independent and identically distributed (i.i.d.) for each type of patients, bed occupancy times are not i.i.d.\ because they depend on bed admission times. In the simulation model, we assume that the lengths of stay of patients of each type follow the empirical distribution obtained from the data, and that the discharge times of all patients follow the distribution specified in Figure~\ref{fig:discharge-hours}.

Aside from the availability of beds, the bed assignment and patient transfer process may also be affected by secondary resources such as BMU staff and nurses in the ED and wards. To faithfully reproduce patient flow dynamics, we add pre- and post-allocation delays to patients' waiting times, following the simulation model by \citet{ShiETAL16}. From the instant when a bed request is received till the point when the patient is assigned to an available bed, several activities can be performed by the BMU staff, such as searching through the bed management system, batching the bed requests for processing later, liaising with ward staff and patients whenever operational needs arise. A patient's pre-allocation delay is the duration of this process. Once the patient is assigned to an available bed, the BMU will initiate the transfer process that entails handing over the patient's records to ward staff and moving the patient from the ED to the ward. A patient's post-allocation delay is the duration of the transfer process. Please refer to Section~3.4 in \citet{ShiETAL16} for more details. We illustrate the mean pre- and post-allocation delays from the records of ED and SOC patients in Figure~\ref{fig:allocation} in Section~\ref{sec:more-simulation}, along with the corresponding standard deviations. In the simulation model, we assume that pre- and post-allocation delays follow log-normal distributions with time-varying means and standard deviations specified in Figure~\ref{fig:allocation}, as \citet{ShiETAL16} did in their paper.

\subsection{The Hospital's Bed Assignment Policy}
\label{sec:hospital-policy}

To replicate performance statistics from the patient flow data, we must faithfully represent the hospital's bed assignment practice in the simulation model. Our partner hospital is relatively conservative toward patient overflowing, always trying their best to minimize overflow rates. The following guidelines are recommended for their inpatient bed management.

When a bed request is received and a primary bed is found available, the BMU should assign the patient to the primary bed immediately. If there is no primary bed, the BMU may either assign the patient to a primary bed to be available later, or put the bed request on hold until a primary bed is available or overflowing is triggered. We adopt the method proposed by \citet{ShiETAL16} to simulate the latter cases: If a bed request is raised at time $ t $ of a day and no primary bed is available, then with probability $ p(t) $, the patient will be assigned to a primary bed to be available later, and with probability $ 1 - p(t) $, the bed request will be put on hold. In simulation, we take 
\[  
p(t) = 
\begin{cases}
0 & \mbox{for $ 0\leq t < 15 $,}\\
0.7 & \mbox{for $ 15 \leq t <21 $,}\\
0 & \mbox{for $ 21 \leq t <24 $}
\end{cases}
\]
for all days, where $ t $ is the time of the day measured in hours. This formula can be roughly verified by the data and interpreted according to the BMU's practice. From 3\,\textsc{pm} to 9\,\textsc{pm}, some patients are assigned to beds to be available later, because most beds vacated on the present day have been taken by waiting patients by 3\,\textsc{pm}. The BMU would thus reserve some beds for incoming patients to control the overflow rate.

When a bed becomes available and a patient has been assigned to this bed, the BMU will initiate the patient transfer process immediately. If no patient has been assigned to this bed and there are primary patients on hold, the BMU will generally follow the FCFS rule, assigning the primary patient who has the earliest bed request time to this bed. If no patient has been assigned to the bed and there is no primary patient waiting, the BMU will leave the bed vacant until a primary patient's bed request is received or patient overflowing is triggered.

The BMU begins to conduct patient overflowing at 5\,\textsc{pm} each day. Bed requests received before 3\,\textsc{pm} and still on hold will be processed in the FCFS manner. If a non-primary bed is available, a patient will be sent there immediately. If the BMU cannot find a bed for a patient, the bed request will be on hold again. Since primary beds are unlikely to be available later, patient overflowing will be triggered repeatedly at the beginning of each hour after 5\,\textsc{pm}, until all patients who requested beds before 3\,\textsc{pm} are assigned to beds. While bed requests received after 3\,\textsc{pm} are eligible for overflowing on the next day, most of them will be assigned to primary beds before the next day's overflowing is triggered.

\subsection{Model Validation}
\label{sec:validation}

\begin{figure}[t]
\centering
\begin{subfigure}{0.49\textwidth}
\includegraphics[trim={.3in 2.5in .4in 2.5in},height=2.35in]{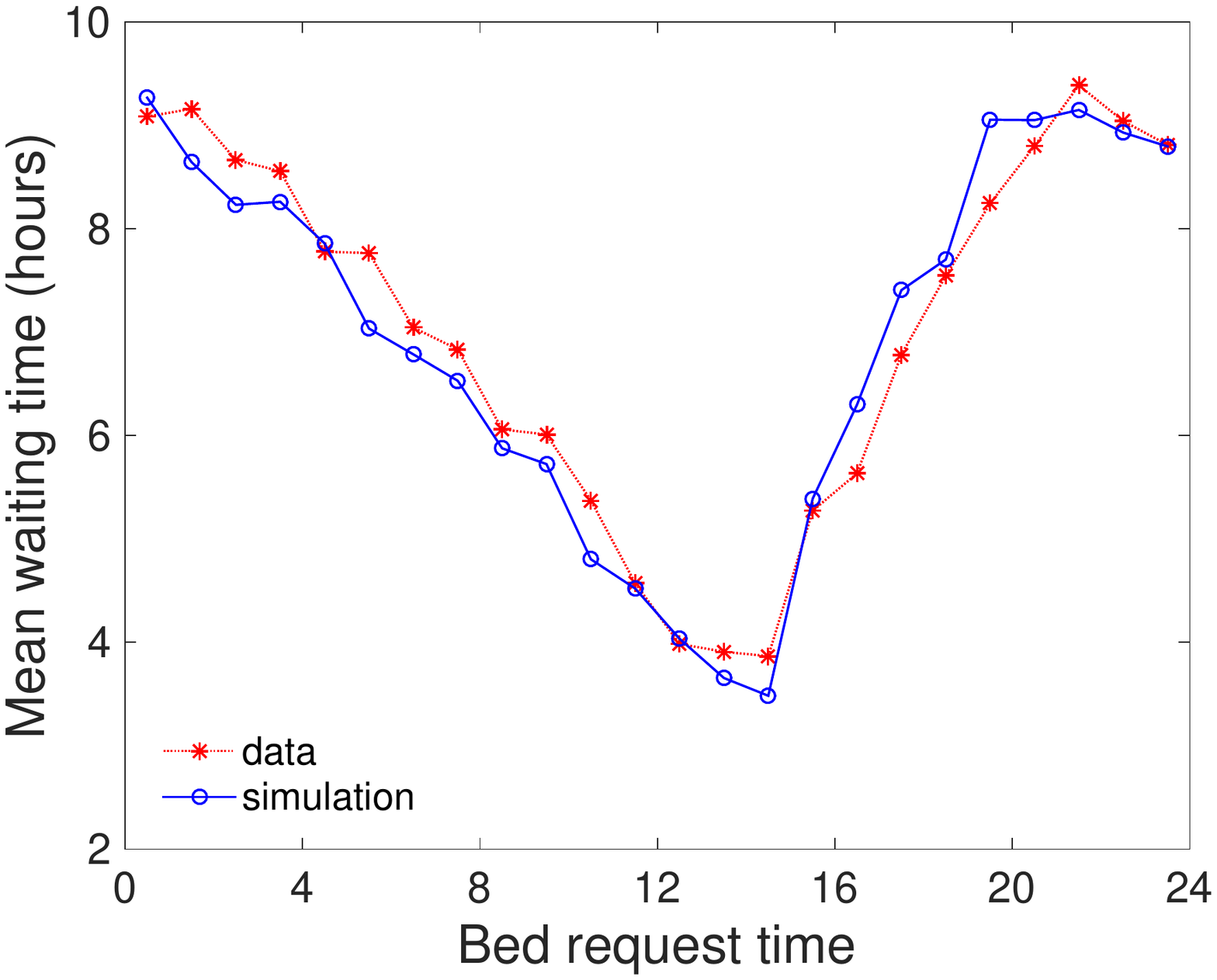}
\caption{Mean waiting times}
\label{fig:validation-mean}
\end{subfigure}
\begin{subfigure}{.49\textwidth}
\includegraphics[trim={.3in 2.5in .4in 2.5in},height=2.35in]{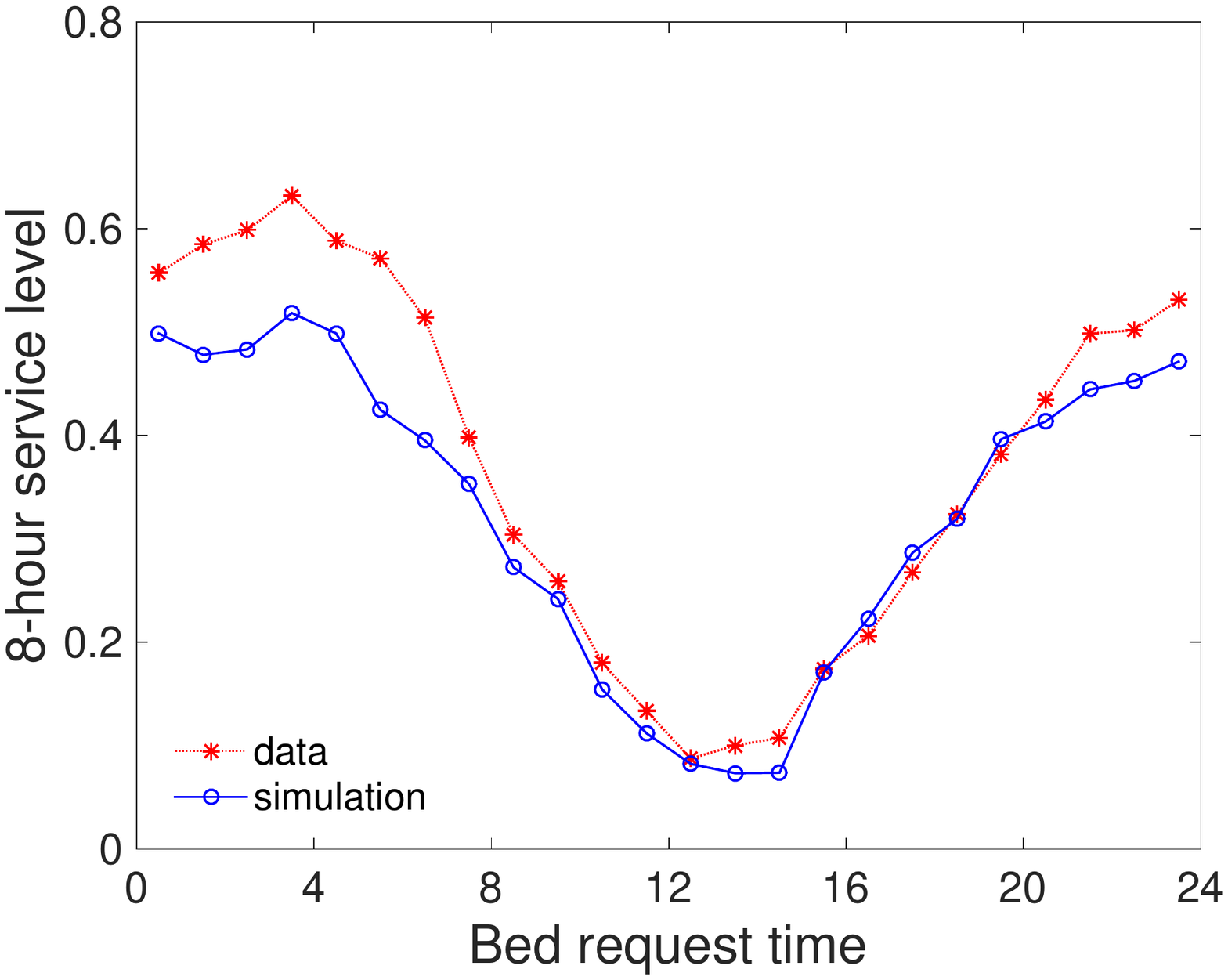}
\caption{8-hour service level}
\label{fig:validation-8}
\end{subfigure}
\begin{subfigure}{.49\textwidth}
\includegraphics[trim={.3in 2.5in .4in 2.5in},height=2.35in]{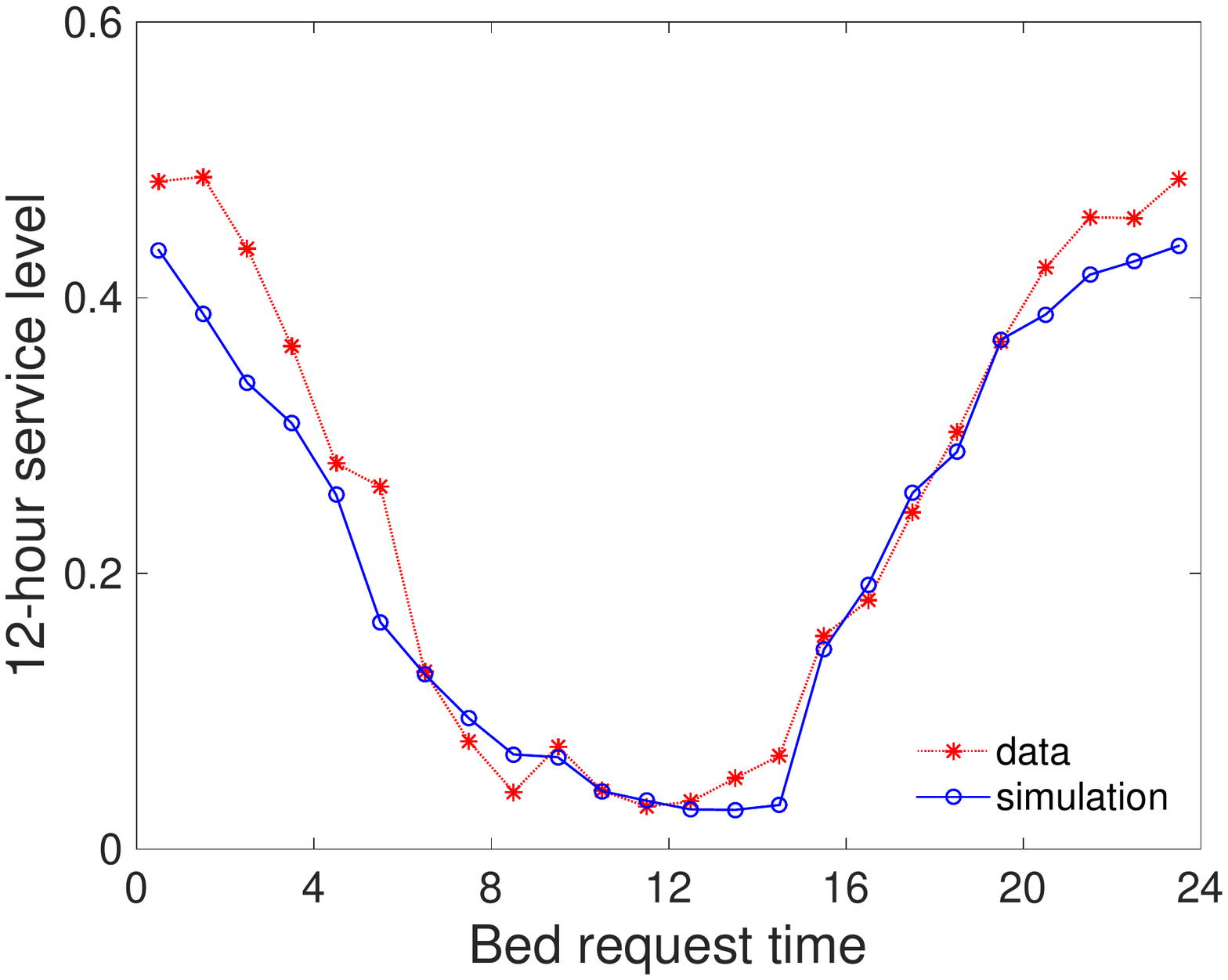}
\caption{12-hour service level}
\label{fig:validation-12}
\end{subfigure}
\begin{subfigure}{.49\textwidth}
\includegraphics[trim={.3in 2.5in .4in 2.5in},height=2.35in]{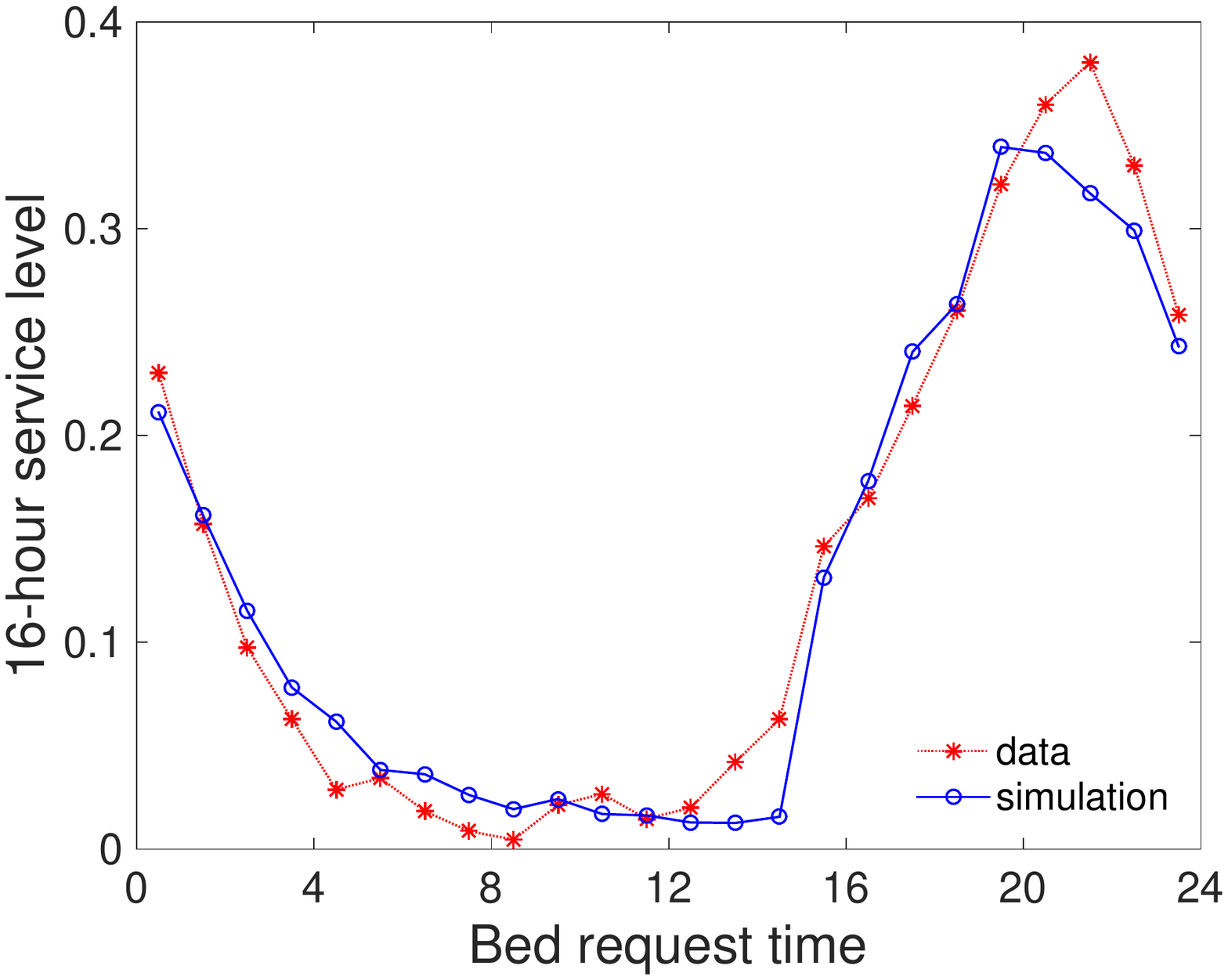}
\caption{16-hour service level}
\label{fig:validation-16}
\end{subfigure}
\medskip
\caption{Performance comparison between statistics from the data and simulation results.}
\label{fig:validation}
\end{figure}

To validate the simulation model, we compare some performance estimates obtained by simulation with corresponding statistics from the patient flow records. Using a software platform coded in Python, we simulate the hospital's inpatient operations for 1000 days, excluding the first 100 days as the burn-in period. The mean waiting times and several service levels of ED patients across the day are plotted in Figure~\ref{fig:validation}. Here, a patient's waiting time refers to the duration from the bed request time to the bed admission time, including both pre- and post-allocation delays, and the $ n $-hour service level is the percentage of patients whose waiting times exceed $ n $ hours. This figure clearly shows that our simulation model can generally capture the hospital's performance in terms of these metrics. The mean waiting times from the ED to various bed pools are reported in Table~\ref{tab:validate-delays} in Section~\ref{sec:more-simulation}, where the simulation results generally agree with the results from the data for large pools (e.g., the pools of B2 and C beds in medicine and surgery).

\begin{figure}[t]
\centering
\begin{minipage}{.49\textwidth}
\centering
\includegraphics[trim={.3in 2.5in .4in 2.5in},height=2.35in]{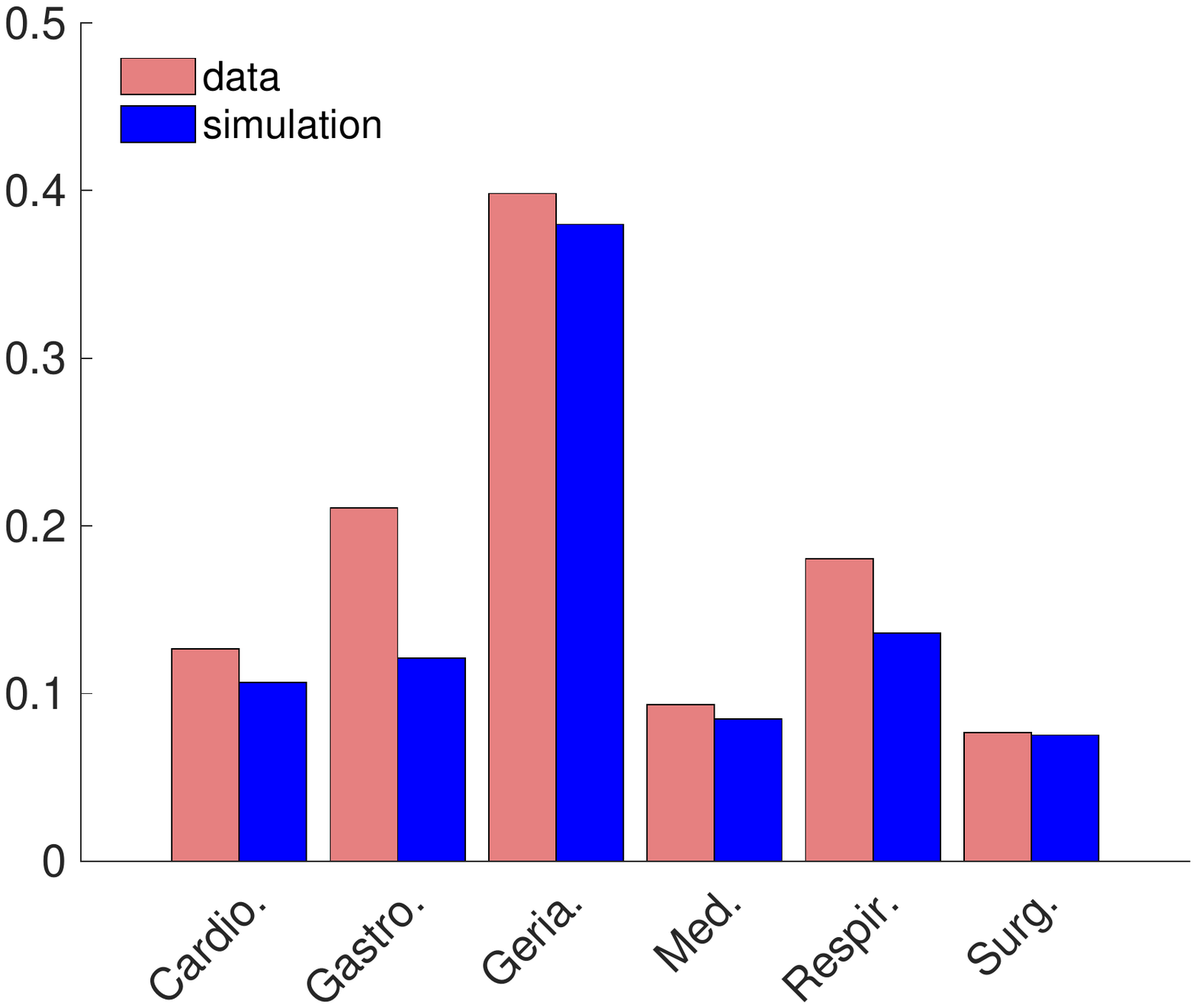}
\caption{Overflow rates in different specialties.}
\label{fig:overflow}
\end{minipage}~~
\begin{minipage}{.49\textwidth}
\centering
\includegraphics[trim={.3in 2.5in .4in 2.5in},height=2.35in]{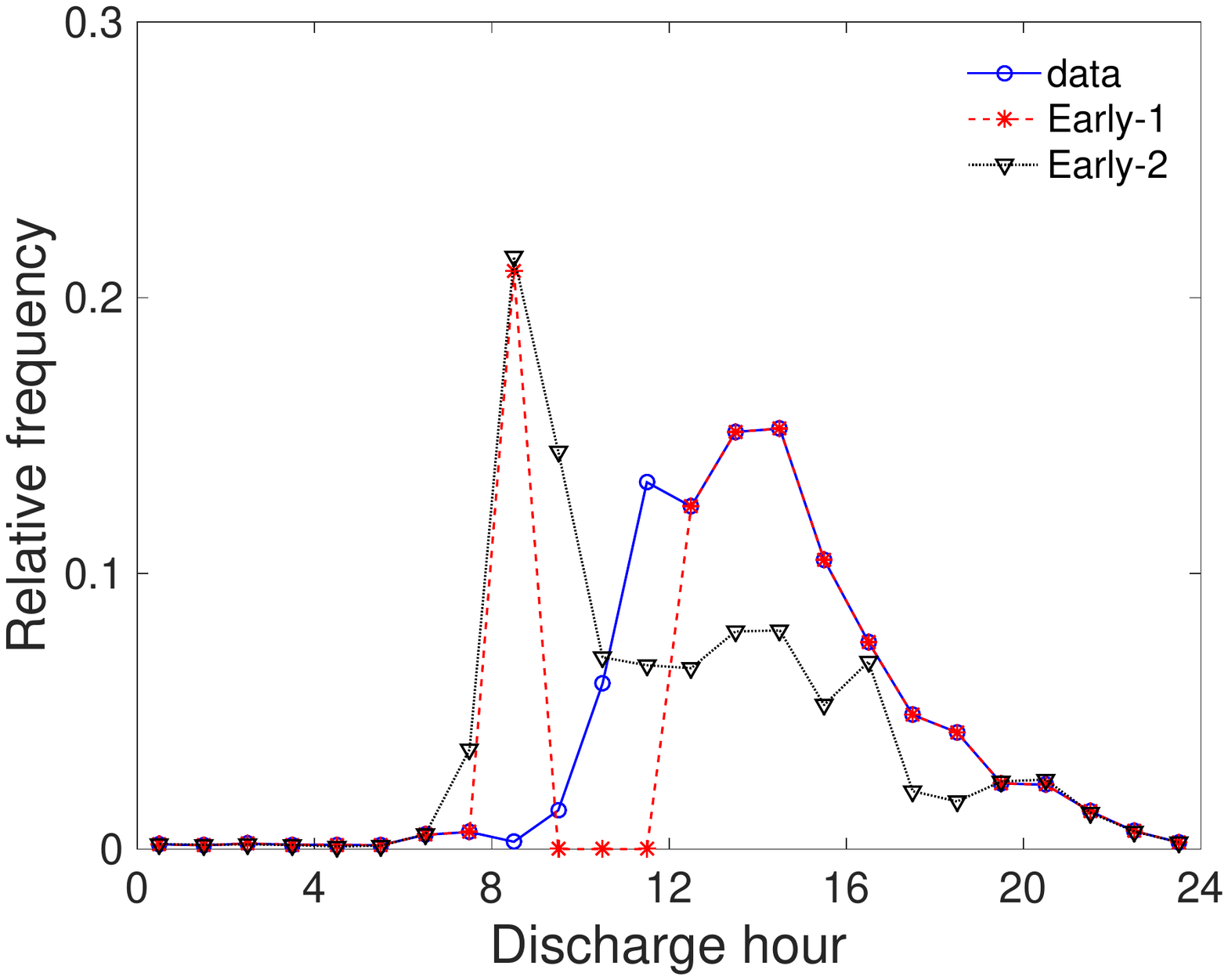}
\caption{Early discharge time distributions.}
\label{fig:early-discharge-distribution}
\end{minipage}
\end{figure}

We compare overflow rates in different specialties in Figure~\ref{fig:overflow}. The overflow rates produced by simulation are generally close to those from the data, except in geriatrics and respiratory medicine. These two specialties are relatively small, having 24 and 43 beds, respectively. We would admit that our simulation model are unable to represent all bed assignment decisions in practice. For example, the BMU staff may occasionally send patients in a serious condition to non-primary beds before 5\,\textsc{pm}, to ensure that the patients can receive inpatient care as soon as possible. Such special cases are not captured in our simulation model. As a result, we may slightly underestimate the overflow rate by simulation. The overflow rate of all ED patients from the data is 12.25\%, while that produced by simulation is 10.88\%. 

\section{Numerical Assessment of the $ \boldsymbol{P} $ Model Approach}
\label{sec:numerical}

We conduct a numerical study in this section to assess the proposed bed assignment approach, using the simulation model introduced in Section~\ref{sec:simulation} as the platform. We consider three performance measures, including the mean waiting times and 12-hour service levels of ED patients and the overflow rates of all patients, for comparison between the proposed approach and several bed assignment approaches commonly used in practice. Although we consider the first two performance measures only for ED patients, we do not assume that these patients have priority over patients from other sources in the bed assignment process. In our simulation, all bed requests are scheduled by the same assignment policy. 

We repeat the simulation process described in Section~\ref{sec:validation} with different bed assignment approaches. In particular, the $ P $ model approach is implemented by calling Gurobi Optimizer 6.5.2. In this approach, a decision iteration is triggered when a bed is available, when a bed request is received and there is at least one available bed for the patient, or when a patient's boarding time is about to exceed the delay target. The overflow cost is taken to be $ u_{ij} = 1$ if bed~$ j $ is a non-primary bed for patient~$ i $, so the overflow budget turns out to be the maximum number of non-primary assignments allowed in an iteration. We first solve \eqref{eq:feasibility} for the minimum overflow budget, next follow \eqref{eq:overflow-budget-integer} to obtain the overflow budget with some $ \alpha $, $ \beta $, and $ \Delta $, and then solve \eqref{eq:optimal-assignment} for an optimal assignment plan. This optimal assignment plan is used to determine the next patient for the available bed. If there is a patient whose boarding time has already exceeded the delay target, we will use the assignment plan obtained by solving \eqref{eq:feasibility} to determine the next patient for the available bed. We include pre- and post-allocation delays in simulation for all bed assignment approaches, in order for the comparison to be fair. Therefore, a patient's waiting time is the sum of the boarding time and the pre- and post-allocation delays. Since we would reduce the 12-hour service level, delay targets for all patients are set to be 10 hours in the $ P $ model approach. We use $ P(\alpha,\beta,\Delta) $ to denote the $ P $ model policy with parameters $ \alpha $, $ \beta $, and $ \Delta $ in \eqref{eq:overflow-budget-integer}. Then for $ \Delta \geq 0 $, $ P(0,0,\Delta) $ refers to the $ P $ model policy using minimum overflow budgets in all iterations.

\subsection{Comparison with the Hospital's Current Practice and Early Discharge Policies}
\label{sec:compare-early}

\begin{figure}[t]
\centering
\begin{subfigure}{0.49\textwidth}
\includegraphics[trim={.3in 2.5in .4in 2.5in},height=2.35in]{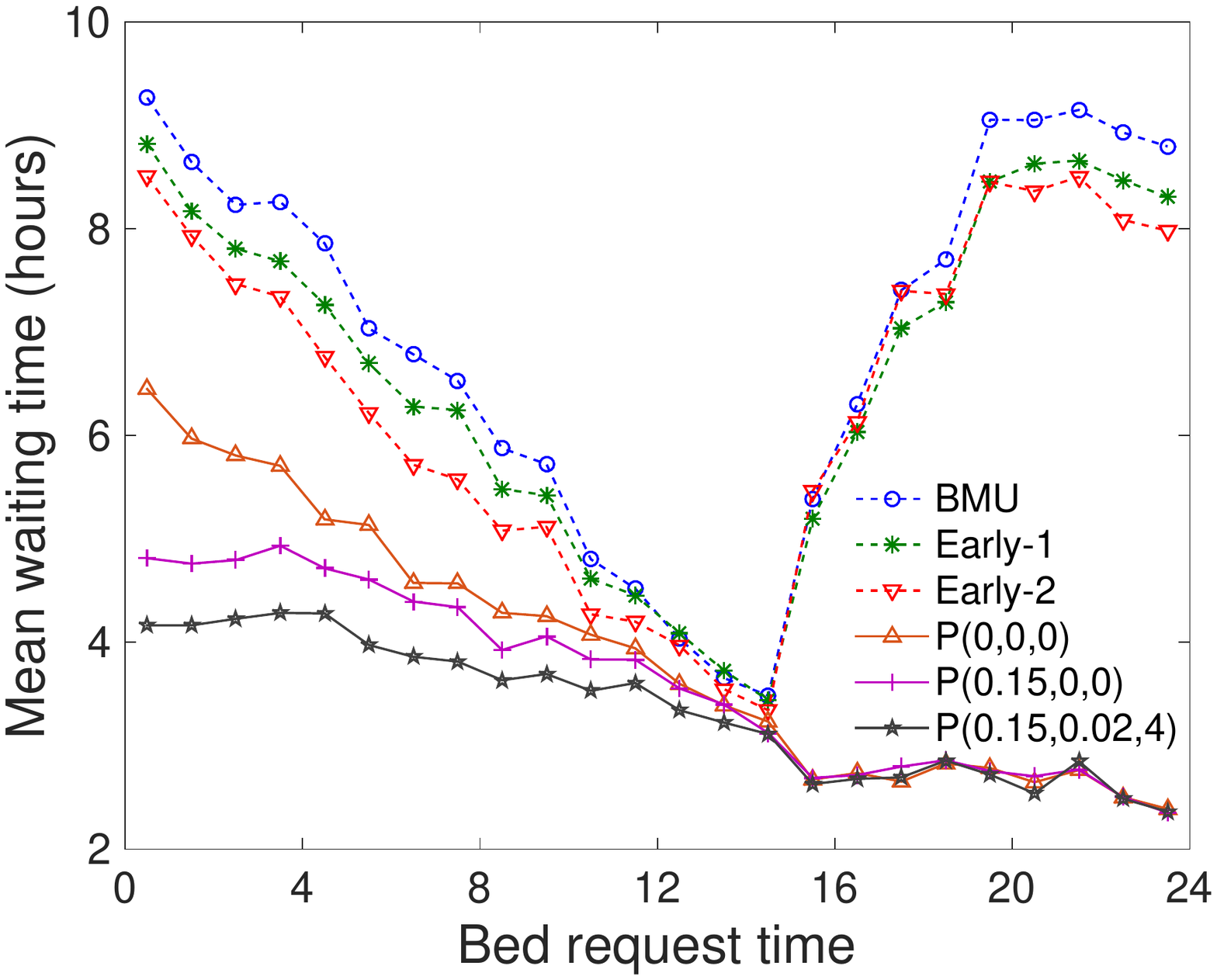}
\caption{Mean waiting times}
\label{fig:means}
\end{subfigure}
\begin{subfigure}{.49\textwidth}
\includegraphics[trim={.3in 2.5in .4in 2.5in},height=2.35in]{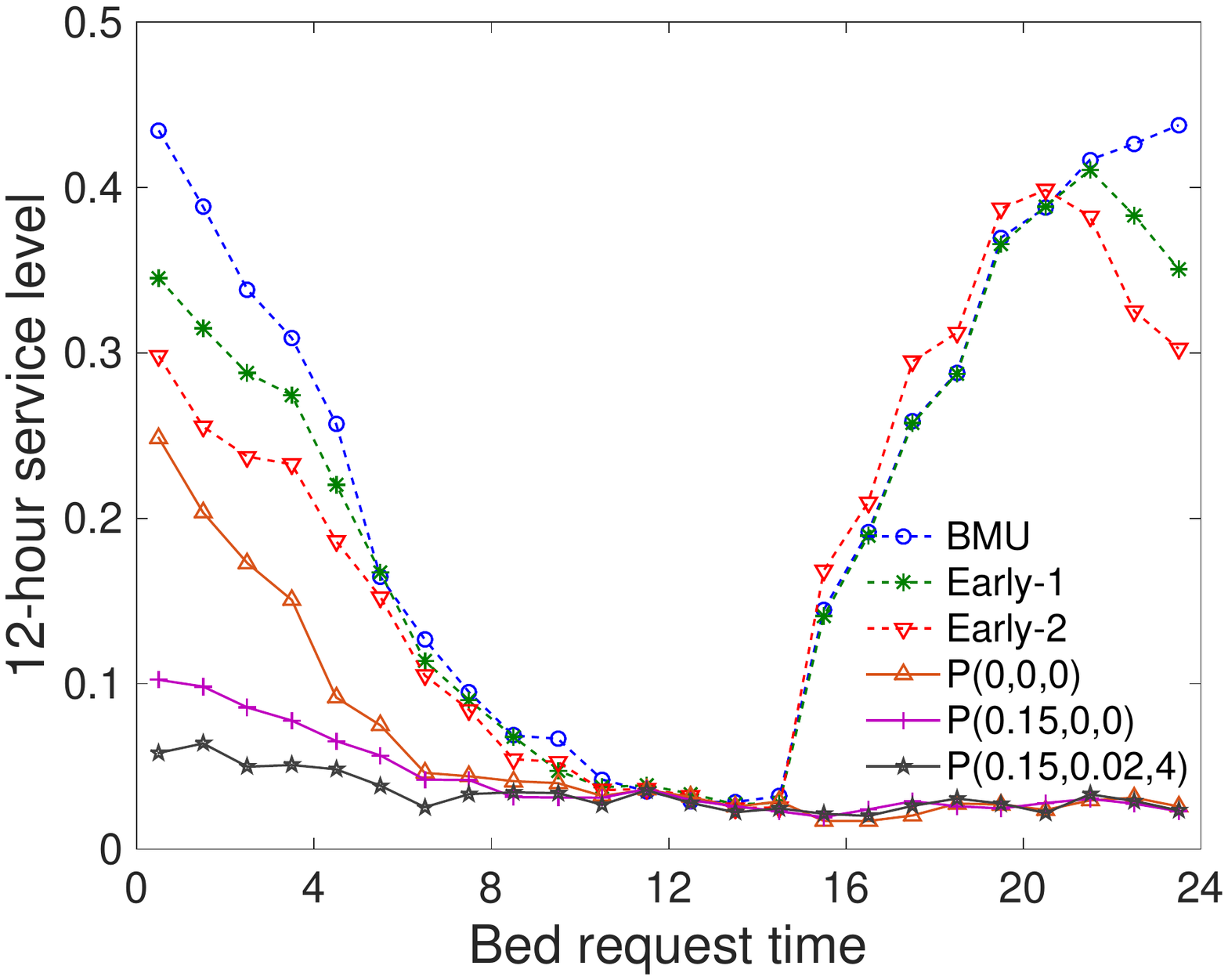}
\caption{12-hour service levels}
\label{fig:service-levels}
\end{subfigure}
\medskip
\caption{Performance comparison between early discharge policies and the $ \boldsymbol{P} $ model approach.}
\label{fig:early-comparison}
\end{figure}

In the first numerical example, we compare the $ P $ model approach with the hospital's current practice as well as two early discharge policies. More specifically, the following bed assignment policies are considered.

\begin{description}
\item[The BMU's policy]: This policy is designed according to the hospital's current guidelines (see Section~\ref{sec:hospital-policy}), denoted by BMU in figures and tables.

\item[The BMU's policy with the first early discharge distribution]: To mitigate the time-of-day effect, \citet{ShiETAL16} proposed a hypothetical policy under which 26\% of patients are discharged before noon and the first discharge peak time appears between 8\,\textsc{am} and 9\,\textsc{am} (see Figure~2b in their paper). They demonstrated that if the distributions of pre- and post-allocation delays are time-invariant, this early discharge policy is capable of equalizing waiting times for ED patients in their simulation model. We consider a similar discharge time distribution, allowing 26\% of patients to be discharged before noon and the first discharge peak time to appear before 9\,\textsc{am} (see Figure~\ref{fig:early-discharge-distribution}). We use Early-1 to denote the BMU's bed assignment policy with this early discharge time distribution.

\item[The BMU's policy with the second early discharge distribution]: We follow the BMU's bed assignment policy, but discharging patients even more aggressively in the morning. We assume that 55\% of patients are discharged before noon and 26\% of patients are discharged before 9\,\textsc{am} (see Figure~\ref{fig:early-discharge-distribution}). This policy is denoted by Early-2 in figures and tables.
\end{description}

We plot the mean waiting times and 12-hour service levels for the above bed assignment policies in Figure~\ref{fig:early-comparison}. The waiting times of ED patients are shortened if a significant portion of patients are discharged in the morning. A more aggressive early discharge policy can further reduce patients' waiting times. However, it turns out that under the BMU's current practice, discharging patients several hours earlier does not suffice to synchronize the patient discharge process with the bed request process. Even if more than one half of patients are discharged before noon, the reduction in the mean waiting time is merely around one hour for patients who request beds at night. Many of them may still wait more than eight hours for beds. As we discussed in Section~\ref{sec:introduction}, it is generally difficult to discharge a significant portion of patients before noon, because physicians and nurses are busy with ward rounds in the morning. For many hospitals, an aggressive early discharge policy is not a viable option.

By maximizing the joint probability of all patients meeting their delay targets in each iteration, the $ P $ model approach can reduce waiting times significantly. Even if the overflow budget is kept minimum in all iterations (i.e., $ \alpha = \beta = 0 $), the mean waiting times can still be reduced up to 6.4 hours at night. As a result, the time-of-day effect on ED boarding is greatly mitigated. If a larger overflow budget is allowed, we may further equalize the waiting times throughout the day. We also test the $ P $ model approach with $ \alpha = 0.15 $ and $ \beta = 0 $. In this case, we may send up to 15\% of waiting patients to non-primary beds when $ B^{\star} $ is small. This adjustment allows us to further reduce waiting times in the early morning at the expense of an increased overflow rate. By setting $ \alpha = 0.15 $, $ \beta = 0.02 $, and $ \Delta = 4 $ hours, we may further improve the delay performance of the $ P $ model approach. In this case, we also adjust the overflow budget according to the bed request rate in the next four hours.

\begin{table}[t]
\centering
\caption{Overflow rates under various bed assignment policies.}
\label{tab:overflow}
\medskip
\begin{tabular} {c c c c c c c c c}
\toprule
& BMU & Early-1 & Early-2 & $ P(0,0,0) $ & $ P(.15,0,0) $ & $ P(.15,.02,4) $ & TB-1 & TB-2 \\ 
\midrule
Overflow rate &  10.89\% & 10.76\% & 10.55\% & 13.12\% & 14.11\% & 15.52\% & 19.50\% & 15.96\%\\
\bottomrule
\end{tabular}
\end{table}

\begin{figure}[t]
\centering
\begin{minipage}{.49\textwidth}
\centering
\includegraphics[trim={.3in 2.5in .4in 2.5in},height=2.35in]{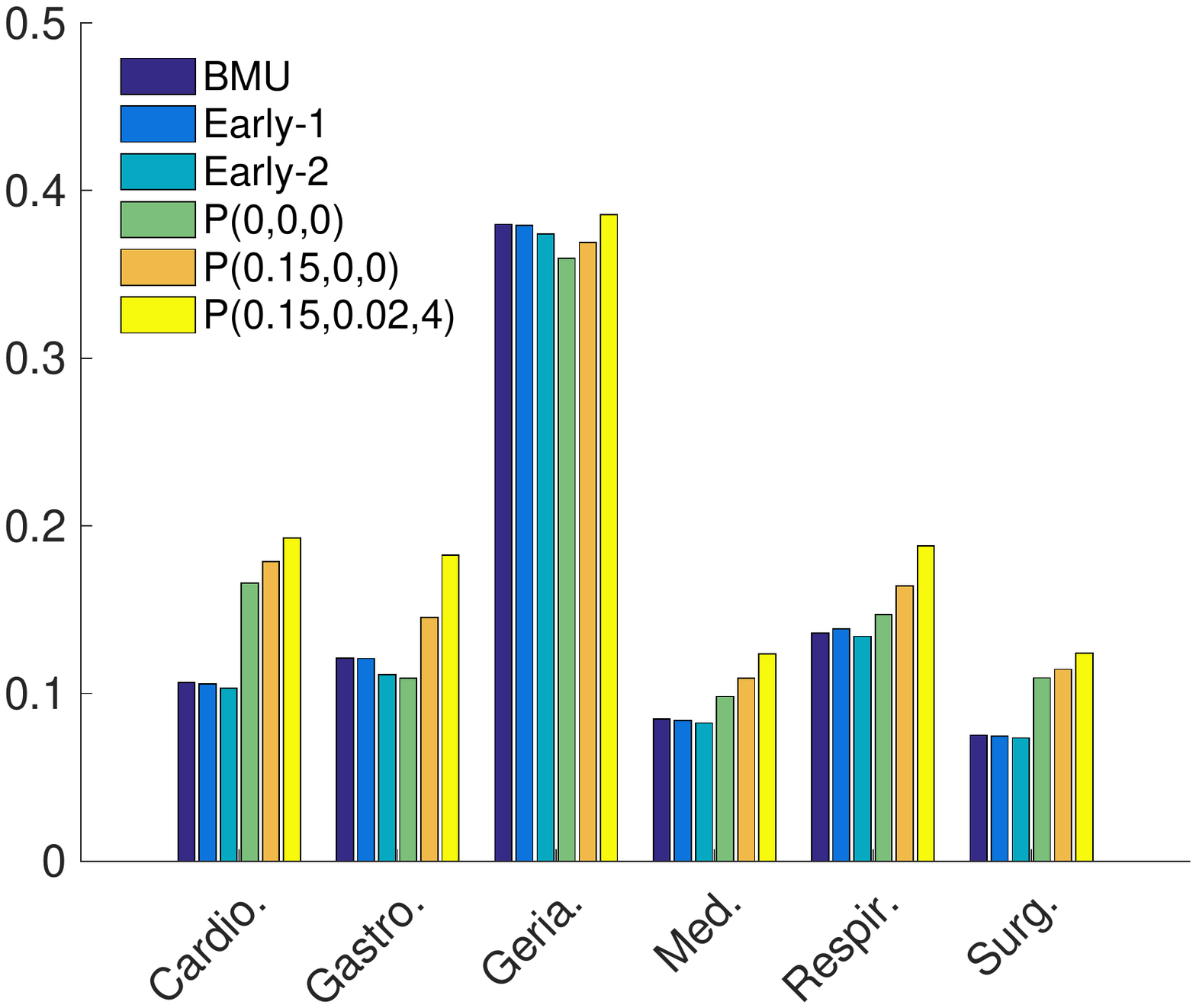}
\caption{Overflow rates in different specialties under various bed assignment policies.}
\label{fig:overflow-comparison}
\end{minipage}~~
\begin{minipage}{.49\textwidth}
\centering
\includegraphics[trim={.3in 2.5in .4in 2.5in},height=2.35in]{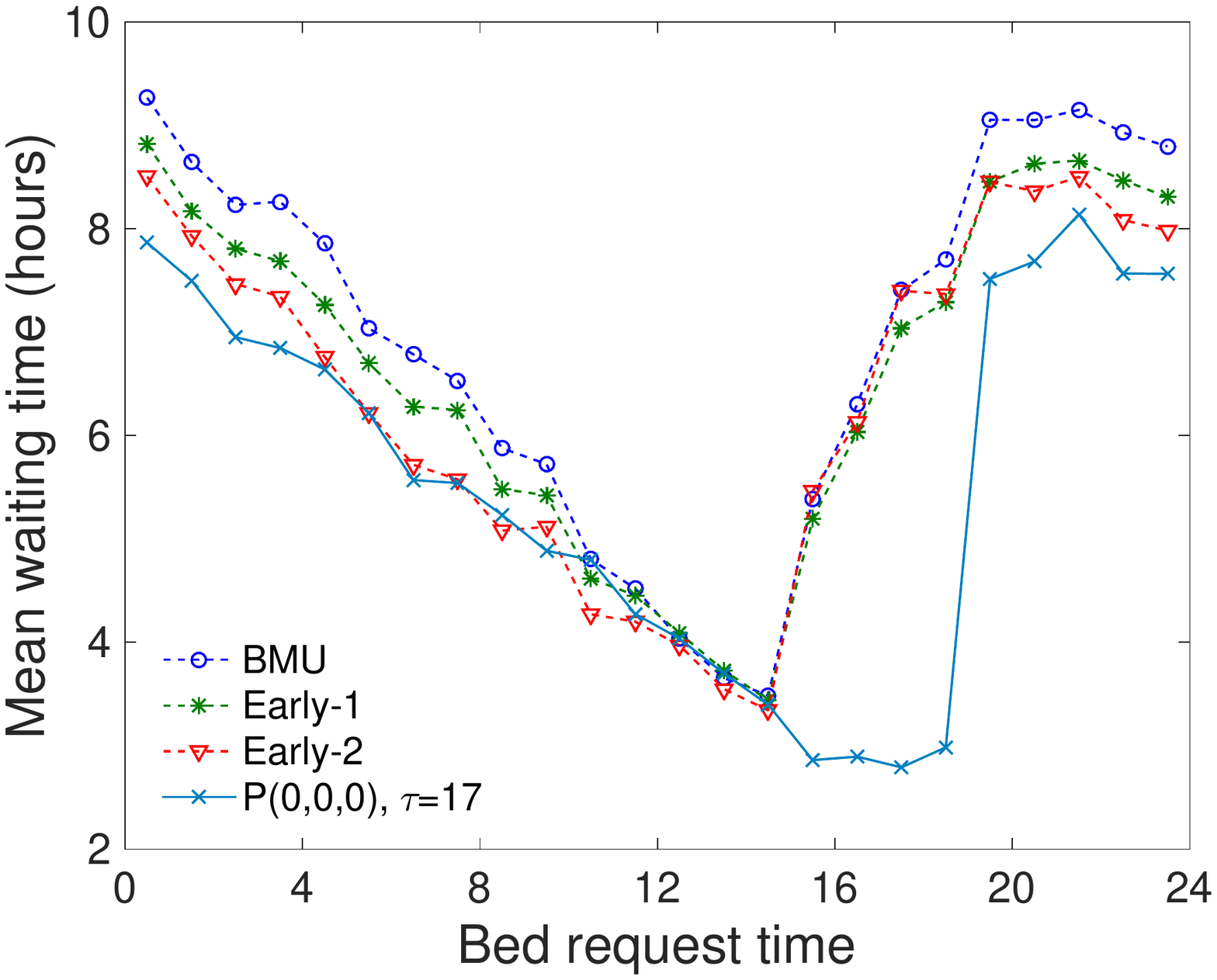}
\caption{The mean waiting times of the $ \boldsymbol{P} $ model formulation with a matched overflow rate.}
\label{fig:matched-overflow}
\end{minipage}
\end{figure}

We report the total overflow rates under these bed assignment policies in Table~\ref{tab:overflow} and illustrate the overflow rates in each specialties in Figure~\ref{fig:overflow-comparison}. The total overflow rate is slightly decreased under the early discharge policies, since more primary beds would be available before patient overflowing is triggered. As we discussed in Section~\ref{sec:hospital-policy}, the hospital's current bed assignment guidelines are generally conservative toward patient overflowing. In this table, the overflow rate under the $ P $ model approach is noticeably larger than that under the BMU's policy, even though the overflow budget takes the minimum values in all iterations. This implies that the BMU's policy would be too conservative if we take 10 hours as the delay target (excluding pre- and post-allocation delays). Since more patients are sent to non-primary beds by the $ P $ model approach, one may raise the question whether it is fair to compare its performance with that of the BMU's policy. To address this concern, we increase the delay target to $ \tau_{i} = 17 $ hours for all patients and take $ \alpha = \beta = 0 $ in \eqref{eq:overflow-budget-integer}. With these parameters, the $ P $ model approach yields an overflow rate of 10.67\%, comparable to the rate of 10.89\% by the BMU's policy. We plot the mean waiting times under this policy in Figure~\ref{fig:matched-overflow}. Although the $ P $ model approach cannot mitigate the time-of-day effect on boarding times with such a large delay target, the improvement in waiting times is apparent compared with the performance of the BMU's policy. For the convenience of illustration, we re-plot the mean waiting times under the two early discharge policies in Figure~\ref{fig:matched-overflow}, where the $ P $ model approach outperforms both of them with a comparable overflow rate. One may further improve the performance of the $ P $ model approach by adjusting the parameters, while maintaining the overflow rate at a similar level. We would omit such an experiment since it is not our focus.

\begin{figure}[t]
\centering
\includegraphics[trim={.3in 2.5in .4in 2.5in},height=2.35in]{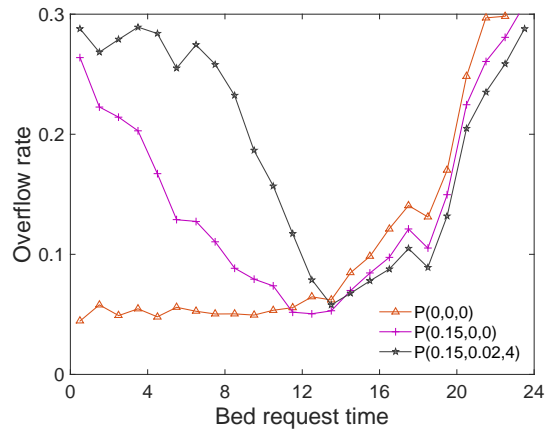}
\caption{Hourly overflow rates under the three $ \boldsymbol{P} $ model policies.}
\label{fig:hourly-overflow}
\end{figure}

We plot hourly overflow rates under the three $ P $ model policies in Figure~\ref{fig:hourly-overflow}. The hourly overflow rate under the $ P $ model approach with $ \alpha = \beta = 0 $ increases with time throughout the day. This is because in our implementation, we assume that the set $ \mathcal{J} $ includes the beds that are either available at the present time or to be available on the present day. At the beginning of a day, all beds to be available on the present day are included in $ \mathcal{J} $. A bed request received in the early morning is thus more likely to be assigned to a primary bed that will be available later. As a result, the overflow rate in the early morning would be relatively low compared with later hours on the same day, while waiting times in the early morning would be relatively long. We may mitigate this issue using larger $ \alpha $ and $ \beta $, which allows us to send more patients to non-primary beds in the morning, thus equalizing boarding times across the day.

\subsection{Comparison with Threshold-Based Overflowing Policies}
\label{sec:threshold-policies}

In many hospitals, patient overflowing is triggered on an individual basis when a patient's boarding time exceeds a threshold. We refer to such a policy as a \emph{threshold-based overflowing policy}. Simple threshold-based policies are widely used in practice (\citealp{ShiETAL16,KilincETAL16}). In the second numerical example, we compare the performance of the $ P $ model approach with the following two threshold-based policies.
\begin{description}
\item[The first threshold-based policy]: Since patients requesting beds at late night or in the early morning are likely to experience excessive delays, we should send them to non-primary beds sooner than others. In the first policy, the threshold time for patients requesting beds from 10\,\textsc{pm} to 3\,\textsc{am} is taken to be two hours, while that for other patients is ten hours. Patients cannot be sent to non-primary beds until overflowing is triggered. When a bed becomes available, the eligible primary or non-primary patient who has the earliest bed request time will be assigned to this bed. We use TB-1 to denote this policy.

\item[The second threshold-based policy]: In the second policy, the threshold time for patients requesting beds from 7\,\textsc{pm} to 12\,\textsc{am} is taken to be two hours, while that for other patients is ten hours. The second policy is denoted by TB-2.
\end{description}

\begin{figure}[t]
\centering
\begin{subfigure}{0.49\textwidth}
\includegraphics[trim={.3in 2.5in .4in 2.5in},height=2.35in]{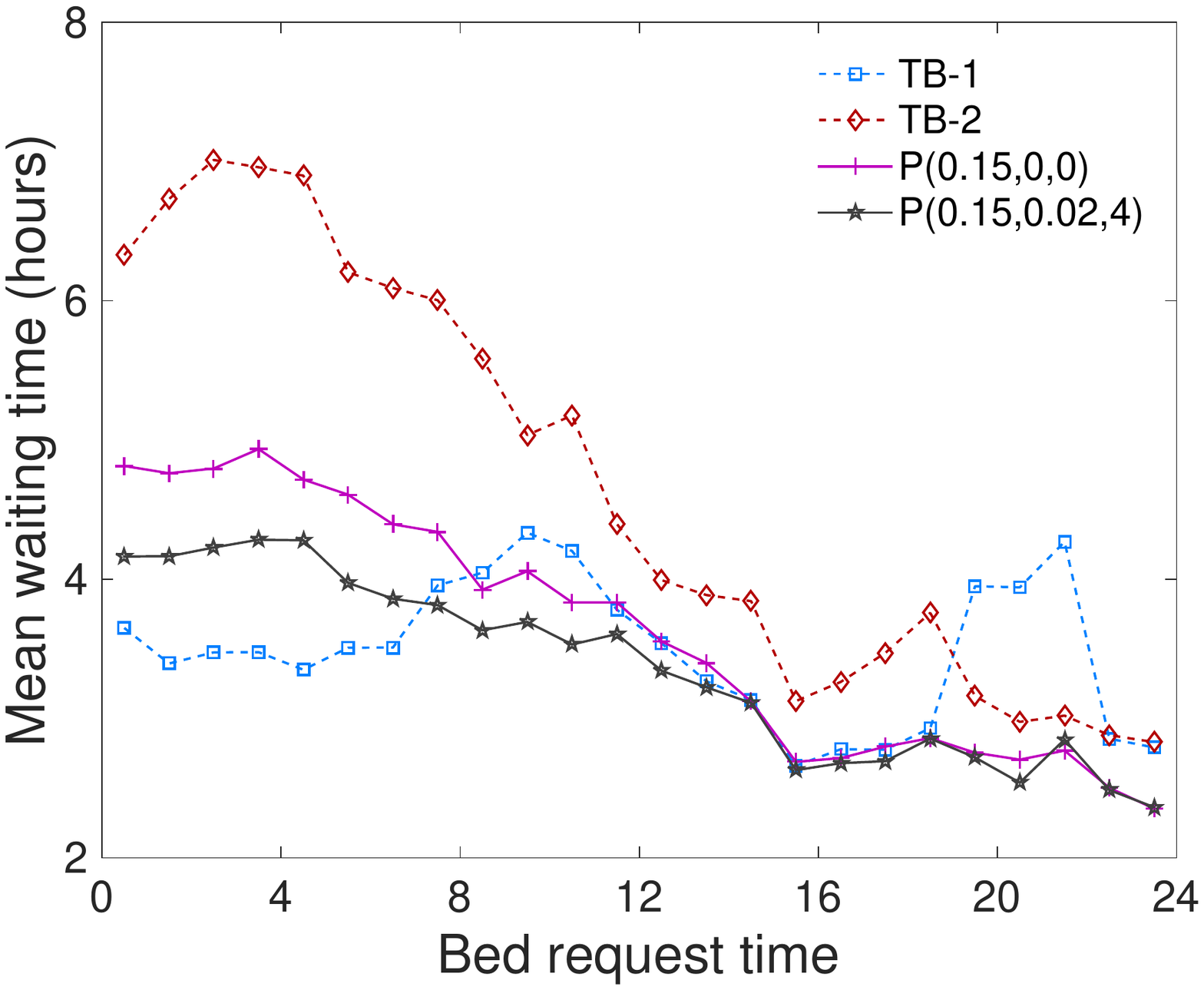}
\caption{Mean waiting times}
\label{fig:delay-means}
\end{subfigure}
\begin{subfigure}{.49\textwidth}
\includegraphics[trim={.3in 2.5in .4in 2.5in},height=2.35in]{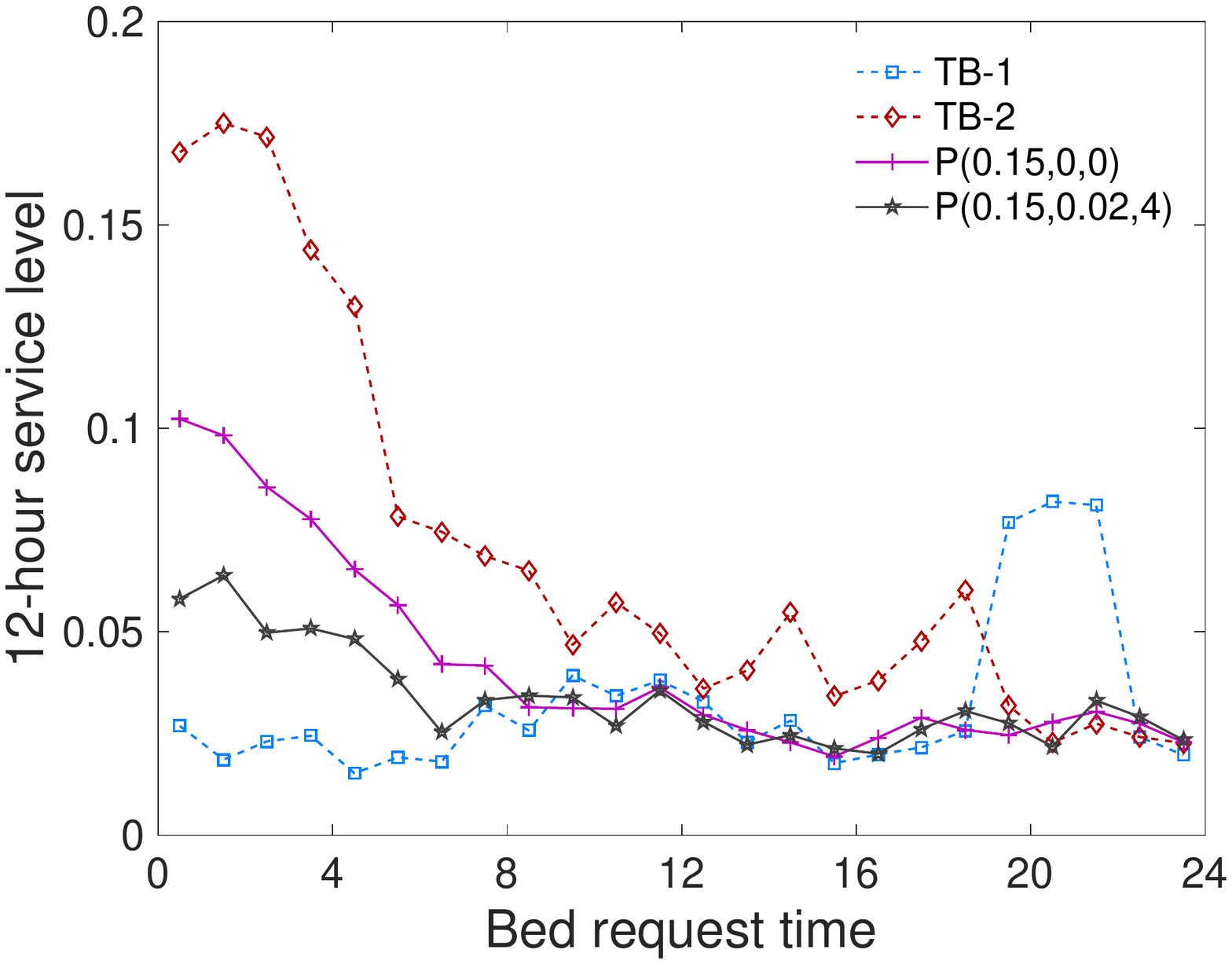}
\caption{12-hour service levels}
\label{fig:delay-service-levels}
\end{subfigure}
\medskip
\caption{Performance comparison between a threshold-based overflowing policy and the $ \boldsymbol{P} $ model approach.}
\label{fig:threshold-policies}
\end{figure}

We plot the mean waiting times and 12-hour service levels for these two threshold-based policies in Figure~\ref{fig:threshold-policies}. The performance of the first policy is comparable to that of the $ P $ model policies. However, the overflow rate reaches 19.50\% under the first threshold-based policy (see Table~\ref{tab:overflow}), much larger than the overflow rates under the $ P $ model policies. We may adjust threshold times to reduce the overflow rate. Under the second threshold-based policy, the overflow rate is 15.96\%, comparable to that under the $ P $ model approach with $ \alpha = 0.15 $, $ \beta = 0.02 $, and $ \Delta = 4 $ hours. However, the waiting time performance under the second policy is not as good as that under the $ P $ model policies. Therefore, by exploiting information about beds to be available later, the $ P $ model approach outperforms threshold-based overflowing policies that are widely used in practice.

\subsection{Computational Performance}
\label{sec:computational}

As we discussed in Section~\ref{sec:P-model}, the $ P $ model approach is computationally tractable because the feasible region of the integer linear program \eqref{eq:optimal-assignment} is sometimes an integral polyhedron. Even if it is not integral, the $ P $ model formulation may still have a ``nice'' structure, allowing us to quickly solve \eqref{eq:optimal-assignment} by a relaxation-based method (such as the cutting plane method) which is embedded in standard integer programming solvers. Although it is beyond our scope to analyze the computational complexity of \eqref{eq:optimal-assignment}, we would demonstrate the computational performance of the $ P $ model approach through numerical experiments.

We repeat the simulation process with the $ P $ model policies and record computation times for solving the $ P $ model formulation in all decision iterations. For the sake of comparison, we also solve the LP relaxation of \eqref{eq:optimal-assignment} in each iteration using the same solver and record the computation time. We report the mean and maximum computation times for solving the $ P $ model and the LP relaxation in Table~\ref{tab:computational}, along with the percentage of iterations in which the LP relaxation has integral solutions. In all three cases, the mean computation time for solving the $ P $ model is just slightly greater than that for solving the LP relaxation, and the maximum computation time for solving the $ P $ model is not significantly greater than that for solving the LP relaxation. If the LP relaxation has an integral solution, it must also be an optimal solution to the $ P $ model. For these three $ P $ model policies, the LP relaxation has integral solutions in around one half of iterations. In general, a larger overflow budget allows the LP relaxation to have integral solutions in more iterations.

\begin{table}[t]
\centering
\caption{Computation times (in seconds) for solving the $ \boldsymbol{P} $ model and the LP relaxation.}
\label{tab:computational}
\medskip
\begin{tabular} {l cc c cc c}
\toprule
& \multicolumn{2}{c}{$ P $ model} & & \multicolumn{3}{c}{LP relaxation} \\
\cmidrule{2-3} \cmidrule{5-7}
& Mean & Max & & Mean & Max & Integral solutions\\\midrule
$ P(0,0,0) $ & 0.034 & 1.500 & & 0.033 & 0.953 & 47.13\%\\
$ P(0.15,0,0) $ & 0.030 & 2.250 & & 0.029 & 1.422 & 52.00\%\\
$ P(0.15,0.02,4) $ & 0.024 & 0.266 & & 0.023 & 0.250 & 59.22\%\\
\bottomrule
\end{tabular}
\end{table} 

\section{Concluding Remarks}
\label{sec:conclusion}

We proposed a data-driven bed assignment approach to balancing boarding and overflowing in hospital wards. The objective is maximizing the percentage of patients whose boarding times are within mandatory targets, without sending excessive patients to non-primary beds. Allowing for critical features of patient flow management, this $ P $ model approach is computationally tractable and easy to implement using existing integer programming solvers. We expect to see more applications of the $ P $ model in the management of large-scale stochastic systems with time-sensitive service requirements.

A data-driven bed assignment platform based on the $ P $ model approach is currently in development by our partner hospital. We may study several more practical issues about inpatient operations using the $ P $ model formulation and incorporate more features into the hospital's platform. For example, how to improve fairness in bed assignment across multiple wards (\citealp{MandelbaumETAL12})? How to schedule elective patients to reduce boarding times in the ED (\citealp{HelmVanOyen14})? How to determine ward partitioning to reduce overflow rates (\citealp{BestETAL15})? We will provide practicable solutions using the $ P $ model approach in the future.

\ACKNOWLEDGMENT{%
The authors would like to thank Melvyn Sim for the constructive discussion that initiated this paper.
}

\bibliographystyle{informs2014} 
\bibliography{refs,healthcare} 

\begin{thebibliography}{35}
\providecommand{\natexlab}[1]{#1}
\providecommand{\url}[1]{\texttt{#1}}
\providecommand{\urlprefix}{URL }

\bibitem[{Armony et~al.(2015)Armony, Israelit, Mandelbaum, Marmor, Tseytlin,
  \protect\BIBand{} Yom-Tov}]{ArmonyETAL15}
Armony M, Israelit S, Mandelbaum A, Marmor YN, Tseytlin Y, Yom-Tov GB (2015) On
  patient flow in hospitals: A data-based queueing-science perspective.
  \emph{Stochastic Systems} 5(1):146--194.

\bibitem[{Armony \protect\BIBand{} Mandelbaum(2011)}]{ArmonyMandelbaum11}
Armony M, Mandelbaum A (2011) Routing and staffing in large-scale service
  systems: The case of homogeneous impatient customers and heterogeneous
  servers. \emph{Operations Research} 59(1):50--65.

\bibitem[{Balinski(1986)}]{Balinski86}
Balinski ML (1986) A competitive (dual) simplex method for the assignment
  problem. \emph{Mathematical Programming} 34(2):125--141.

\bibitem[{Baron et~al.(2017)Baron, Berman, Krass, \protect\BIBand{}
  Wang}]{BaronETAL17}
Baron O, Berman O, Krass D, Wang J (2017) Strategic idleness and dynamic
  scheduling in an open-shop service network: Case study and analysis.
  \emph{Manufacturing \& Service Operations Management} 19(1):52--71.

\bibitem[{Bassamboo et~al.(2006)Bassamboo, Harrison, \protect\BIBand{}
  Zeevi}]{BassambooETAL06}
Bassamboo A, Harrison JM, Zeevi A (2006) Design and control of a large call
  center: Asymptotic analysis of an {LP}-based method. \emph{Operations
  Research} 54(3):419--435.

\bibitem[{Bassamboo \protect\BIBand{} Zeevi(2009)}]{BassambooZeevi09}
Bassamboo A, Zeevi A (2009) On a data-driven method for staffing large call
  centers. \emph{Operations Research} 57(3):714--726.

\bibitem[{Best et~al.(2015)Best, Sandıkçı, Eisenstein, \protect\BIBand{}
  Meltzer}]{BestETAL15}
Best TJ, Sandıkçı B, Eisenstein DD, Meltzer DO (2015) Managing hospital
  inpatient bed capacity through partitioning care into focused wings.
  \emph{Manufacturing \& Service Operations Management} 17(2):157--176.

\bibitem[{Borst et~al.(2004)Borst, Mandelbaum, \protect\BIBand{}
  Reiman}]{BorstETAL04}
Borst S, Mandelbaum A, Reiman MI (2004) Dimensioning large call centers.
  \emph{Operations Research} 52(1):17--34.

\bibitem[{Chan et~al.(2017)Chan, Dong, \protect\BIBand{} Green}]{CDG17}
Chan CW, Dong J, Green LV (2017) Queues with time-varying arrivals and
  inspections with applications to hospital discharge policies.
  \emph{Operations Research} 65(2):469--495.

\bibitem[{Charnes \protect\BIBand{} Cooper(1963)}]{CharnesCooper63}
Charnes A, Cooper WW (1963) Deterministic equivalents for optimizing and
  satisficing under chance constraints. \emph{Operations Research}
  11(1):18--39.

\bibitem[{Dai \protect\BIBand{} Shi(2017{\natexlab{a}})}]{DaiShi18}
Dai JG, Shi P (2017{\natexlab{a}}) Inpatient bed overflow: An approximate
  dynamic programming approach, to appear.

\bibitem[{Dai \protect\BIBand{} Shi(2017{\natexlab{b}})}]{DaiShi17}
Dai JG, Shi P (2017{\natexlab{b}}) A two-time-scale approach to time-varying
  queues in hospital inpatient flow management. \emph{Operations Research}
  65(2):514--536.

\bibitem[{Dai \protect\BIBand{} Tezcan(2008)}]{DaiTezcan08}
Dai JG, Tezcan T (2008) Optimal control of parallel server systems with many
  servers in heavy traffic. \emph{Queueing Systems} 59(2):95--134.

\bibitem[{Green(2008)}]{Green08}
Green LV (2008) Using operations research to reduce delays for healthcare. Chen
  ZL, Raghavan S, eds., \emph{TutORials in Operations Research:
  State-of-the-Art Decision-Making Tools in the Information-Intensive Age},
  1--16 (Hanover, MD: INFORMS).

\bibitem[{Gurvich \protect\BIBand{}
  Whitt(2009{\natexlab{a}})}]{GurvichWhitt09a}
Gurvich I, Whitt W (2009{\natexlab{a}}) Queue-and-idleness-ratio controls in
  many-server service systems. \emph{Mathematics of Operations Research}
  34(2):363--396.

\bibitem[{Gurvich \protect\BIBand{}
  Whitt(2009{\natexlab{b}})}]{GurvichWhitt09b}
Gurvich I, Whitt W (2009{\natexlab{b}}) Scheduling flexible servers with convex
  delay costs in many-server service systems. \emph{Manufacturing \& Service
  Operations Management} 11(2):237--253.

\bibitem[{Gurvich \protect\BIBand{} Whitt(2010)}]{GurvichWhitt10}
Gurvich I, Whitt W (2010) Service-level differentiation in many-server service
  systems via queue-ratio routing. \emph{Operations Research} 58(2):316--328.

\bibitem[{Halfin \protect\BIBand{} Whitt(1981)}]{HalfinWhitt81}
Halfin S, Whitt W (1981) Heavy-traffic limits for queues with many exponential
  servers. \emph{Operations Research} 29(3):567--588.

\bibitem[{Harrison et~al.(2005)Harrison, Shafer, \protect\BIBand{}
  Mackay}]{HarrisonETAL05}
Harrison GW, Shafer A, Mackay M (2005) Modelling variability in hospital bed
  occupancy. \emph{Health Care Management Science} 8(4):325--334.

\bibitem[{He et~al.(2018)He, Sim, \protect\BIBand{} Zhang}]{HeETAL18}
He S, Sim M, Zhang M (2018) Data-driven patient scheduling in emergency
  departments: A hybrid robust--stochastic approach. \emph{Management Science}
  To appear.

\bibitem[{Helm \protect\BIBand{} Van~Oyen(2014)}]{HelmVanOyen14}
Helm JE, Van~Oyen MP (2014) Design and optimization methods for elective
  hospital admissions. \emph{Operations Research} 62(6):1265--1282.

\bibitem[{Huang et~al.(2015)Huang, Carmeli, \protect\BIBand{}
  Mandelbaum}]{HuangETAL15}
Huang J, Carmeli B, Mandelbaum A (2015) Control of patient flow in emergency
  departments, or multiclass queues with deadlines and feedback.
  \emph{Operations Research} 63(4):892--909.

\bibitem[{Kilinc et~al.(2016)Kilinc, Saghafian, \protect\BIBand{}
  Traub}]{KilincETAL16}
Kilinc D, Saghafian S, Traub SJ (2016) Dynamic assignment of patients to
  primary and secondary inpatient units: Is patience a virtue?, preprint.

\bibitem[{Mandelbaum et~al.(2012)Mandelbaum, Mom\v{c}ilovi\'{c},
  \protect\BIBand{} Tseytlin}]{MandelbaumETAL12}
Mandelbaum A, Mom\v{c}ilovi\'{c} P, Tseytlin Y (2012) On fair routing from
  emergency departments to hospital wards: {QED} queues with heterogeneous
  servers. \emph{Management Science} 58(7):1273--1291.

\bibitem[{Nemhauser \protect\BIBand{} Wolsey(1988)}]{NemhauserWolsey88}
Nemhauser GL, Wolsey LA (1988) \emph{Integer and Combinatorial Optimization}
  (New York: Wiley).

\bibitem[{Pines et~al.(2011)Pines, Hilton, Weber, Alkemade, Al~Shabanah,
  Anderson, Bernhard, Bertini, Gries, Ferrandiz, Kumar, Harjola, Hogan, Madsen,
  Mason, \:{O}hl\'{e}n, Rainer, Rathlev, Revue, Richardson, Sattarian,
  \protect\BIBand{} Schull}]{PinesETAL11b}
Pines JM, Hilton JA, Weber EJ, Alkemade AJ, Al~Shabanah H, Anderson PD,
  Bernhard M, Bertini A, Gries A, Ferrandiz S, Kumar VA, Harjola VP, Hogan B,
  Madsen B, Mason S, \:{O}hl\'{e}n G, Rainer T, Rathlev N, Revue E, Richardson
  D, Sattarian M, Schull MJ (2011) International perspectives on emergency
  department crowding. \emph{Academic Emergency Medicine} 18(12):1358--1370.

\bibitem[{Samiedaluie et~al.(2017)Samiedaluie, Kucukyazici, Verter,
  \protect\BIBand{} Zhang}]{SamiedaluieETAL17}
Samiedaluie S, Kucukyazici B, Verter V, Zhang D (2017) Managing patient
  admissions in a neurology ward. \emph{Operations Research} 65(3):635--656.

\bibitem[{Shi et~al.(2016)Shi, Chou, Dai, Ding, \protect\BIBand{}
  Sim}]{ShiETAL16}
Shi P, Chou MC, Dai JG, Ding D, Sim J (2016) Models and insights for hospital
  inpatient operations: Time-dependent {ED} boarding time. \emph{Management
  Science} 62(1):1--28.

\bibitem[{Singer et~al.(2011)Singer, Thode~Jr, Viccellio, \protect\BIBand{}
  Pines}]{SingerETAL11}
Singer AJ, Thode~Jr HC, Viccellio P, Pines JM (2011) The association between
  length of emergency department boarding and mortality. \emph{Academic
  Emergency Medicine} 18(12):1324--1329.

\bibitem[{Stolyar \protect\BIBand{} Tezcan(2011)}]{StolyarTezcan11}
Stolyar AL, Tezcan T (2011) Shadow-routing based control of flexible
  multiserver pools in overload. \emph{Operations Research} 59(6):1427--1444.

\bibitem[{Sun et~al.(2013)Sun, Hsia, Weiss, Zingmond, Liang, Han, McCreath,
  \protect\BIBand{} Asch}]{SunETAL13}
Sun BC, Hsia RY, Weiss RE, Zingmond D, Liang LJ, Han W, McCreath H, Asch SM
  (2013) Effect of emergency department crowding on outcomes of admitted
  patients. \emph{Annals of Emergency Medicine} 61(6):605--611.

\bibitem[{Teow et~al.(2011)Teow, El-Darzi, Foo, Jin, \protect\BIBand{}
  Sim}]{TeowETAL11}
Teow KL, El-Darzi E, Foo C, Jin X, Sim J (2011) Intelligent analysis of acute
  bed overflow in a tertiary hospital in {S}ingapore. \emph{Journal of Medical
  Systems} 36(3):1873--1882.

\bibitem[{Thompson et~al.(2009)Thompson, Nunez, Garfinkel, \protect\BIBand{}
  Dean}]{ThompsonETAL09}
Thompson S, Nunez M, Garfinkel R, Dean M (2009) Efficient short-term allocation
  and reallocation of patients to floors of a hospital during demand surges.
  \emph{Operations Research} 57(2):261--273.

\bibitem[{Wertheimer et~al.(2014)Wertheimer, Jacobs, Bailey, Holstein,
  Chatfield, Ohta, Horrocks, \protect\BIBand{} Hochman}]{WertheimerETAL14}
Wertheimer B, Jacobs REA, Bailey M, Holstein S, Chatfield S, Ohta B, Horrocks
  A, Hochman K (2014) Discharge before noon: An achievable hospital goal.
  \emph{Journal of Hospital Medicine} 9(4):210--214.

\bibitem[{Wertheimer et~al.(2015)Wertheimer, Jacobs, Iturrate, Bailey,
  \protect\BIBand{} Hochman}]{WertheimerETAL15}
Wertheimer B, Jacobs REA, Iturrate E, Bailey M, Hochman K (2015) Discharge
  before noon: Effect on throughput and sustainability. \emph{Journal of
  Hospital Medicine} 10(10):664--669.

\end{thebibliography}

\newpage

%
%
%




\ECSwitch


\ECHead{Appendix}

\section{Proof of Proposition~\ref{prop:FCFS}}
\label{sec:proofs}

\begin{proof}{Proof of Proposition~\ref{prop:FCFS}.}
Let $ t $ be the present time. Then, $ \tilde{d}_{j} = t $. Since $ B \geq B^{\star} $, the feasible region of \eqref{eq:optimal-assignment} is nonempty and an optimal solution exits. Suppose that there is an optimal solution $ \hat{\boldsymbol{z}} $ to \eqref{eq:optimal-assignment} under which patient~$ k $ is assigned to bed~$ j $ but does not have the earliest bed request time among the patients of the same type. In other words, with $ \hat{z}_{kj} = 1$, there exists $ \ell \in \mathcal{I} $ such that $ a_{\ell} < a_{k} $ where
\[ 
a_{\ell} = \min \{a_{i} : i \in \mathcal{I},\ \mathcal{J}_{i} = \mathcal{J}_{k},\ \tau_{i} = \tau_{k},\ u_{im} = u_{k m} \mbox{ for }  m \in \mathcal{J}_{k} \}.
\] 

Let $ h \in \mathcal{J} $ be the bed to which patient~$ \ell $ is assigned under $ \hat{\boldsymbol{z}} $, i.e., $ \hat{z}_{\ell h} = 1 $. Then, $ \tilde{d}_{h} \geq t $ since bed~$ h $ may or may not be available at time $ t $. We modify $ \hat{\boldsymbol{z}} $ by assigning patient~$ k $ to bed~$ h $ and patient~$ \ell $ to bed~$ j $. Let $ \boldsymbol{z}^{\star} $ be the resulting assignment plan. Then,
\[  
z^{\star}_{im} = 
\begin{cases}
1 & \mbox{if $ (i,m) = (k,h) $ or $ (i,m) = (\ell, j) $,}\\
\hat{z}_{im} & \mbox{otherwise.}
\end{cases}
\]
Clearly, $ \boldsymbol{z}^{\star} $ satisfies all constraints in \eqref{eq:optimal-assignment}. Consider the objective function in \eqref{eq:optimal-assignment}, given by
\[  
V(\boldsymbol{z}) = \sum_{i\in\mathcal{I}} \sum_{m\in\mathcal{J}_{i}} z_{im}\ln \mathbb{P}\big(\tilde{d}_{m} - a_{i} \leq \tau_{i}\big).
\] 
Note that $ \tilde{d}_{j} = t \leq \tilde{d}_{h} $. If $ t - a_{\ell} > \tau_{\ell} $, we have $ \mathbb{P}\big(\tilde{d}_{j} - a_{\ell} \leq \tau_{\ell}\big) =  \mathbb{P}\big(\tilde{d}_{h} - a_{\ell} \leq \tau_{\ell}\big) = 0 $. In this case, $ V(\boldsymbol{z}^{\star}) = V(\hat{\boldsymbol{z}}) = -\infty $. If $ t - a_{\ell} \leq \tau_{\ell} $ and $ V(\hat{\boldsymbol{z}}) > -\infty $, we have
\begin{align*}
V(\boldsymbol{z}^{\star}) - V(\hat{\boldsymbol{z}}) & = \ln \mathbb{P}\big(\tilde{d}_{j} - a_{\ell} \leq \tau_{\ell}\big) + \ln \mathbb{P}\big(\tilde{d}_{h} - a_{k} \leq \tau_{k}\big) -
\ln \mathbb{P}\big(\tilde{d}_{j} - a_{k} \leq \tau_{k}\big) - \ln \mathbb{P}\big(\tilde{d}_{h} - a_{\ell} \leq \tau_{\ell}\big) \\
& = \ln \mathbb{P}\big(\tilde{d}_{h} - a_{k} \leq \tau_{k}\big) - \ln \mathbb{P}\big(\tilde{d}_{h} - a_{\ell} \leq \tau_{\ell}\big) \\
& \geq 0.
\end{align*}
Since $ \hat{\boldsymbol{z}} $ is an optimal solution to \eqref{eq:optimal-assignment}, we must have $ V(\boldsymbol{z}^{\star}) = V(\hat{\boldsymbol{z}}) $ and $ \boldsymbol{z}^{\star} $ is also optimal. This implies that by switching the beds for patients~$ k $ and~$ \ell $, we can obtain an optimal assignment plan $ \boldsymbol{z}^{\star} $ that satisfies the condition specified in the proposition. \hfill\Halmos
\end{proof}

\section{More on the Simulation Model}
\label{sec:more-simulation}

More details about the simulation model is provided in this appendix. 

In Table~\ref{tab:LOS}, we report the mean lengths of stay and the corresponding standard deviations of the six specialties from different patient sources. From this table, we can tell that in all six specialties, patients admitted in the morning have shorter mean lengths of stay than those admitted after 12\,\textsc{pm}.

We report the distributions of bed requests for the four accommodation classes in the six specialties in Table~\ref{tab:accommodation}. Beds of classes B2 and C have much higher demand than the other two classes because they are subsidized by the government.

The gender distributions of patients in the six specialties are reported in Table~\ref{tab:gender}. There are more male patients than female patients in all specialties except for geriatrics.

We list the indexes of the 34 bed pools in Table~\ref{tab:bed-pools}, along with the numbers of beds belonging to these pools. In the surgical wards of this hospital, some B2 beds are reserved for EL and SDA patients. These beds are denoted by ``Surg. (EL)'' in Table~\ref{tab:bed-pools}. The hospital also reserves some C beds as two common (secondary) pools for male and female patients in all specialties requesting B2 and C beds. We create a new ``mixed'' specialty in Table~\ref{tab:bed-pools} for these beds.

We list the primary and non-primary bed pools for male and female patients in Tables~\ref{tab:overflow-male} and~\ref{tab:overflow-female}, respectively, according to the hospital's guidelines. The BMU will try to send patients to preferred beds when primary beds are not available. Only when a preferred bed is not available, may the BMU send the patient to a secondary bed.

\begin{figure}[tbh]
\centering
\begin{subfigure}{0.49\textwidth}
\includegraphics[trim={.3in 2.5in .4in 2.5in},height=2.35in]{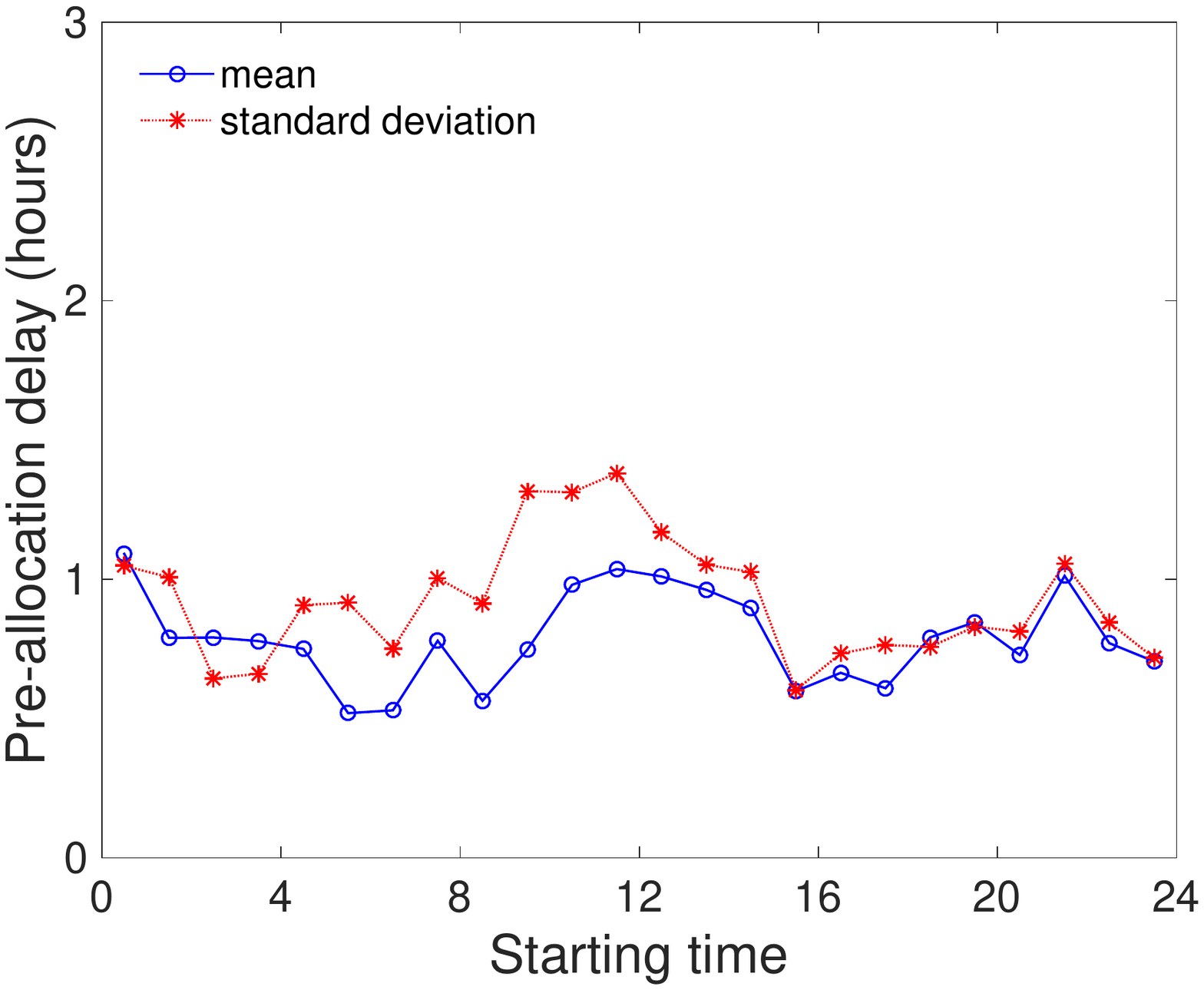}
\caption{Pre-allocation delays}
\label{fig:pre-allocation}
\end{subfigure}
\begin{subfigure}{.49\textwidth}
\includegraphics[trim={.3in 2.5in .4in 2.5in},height=2.35in]{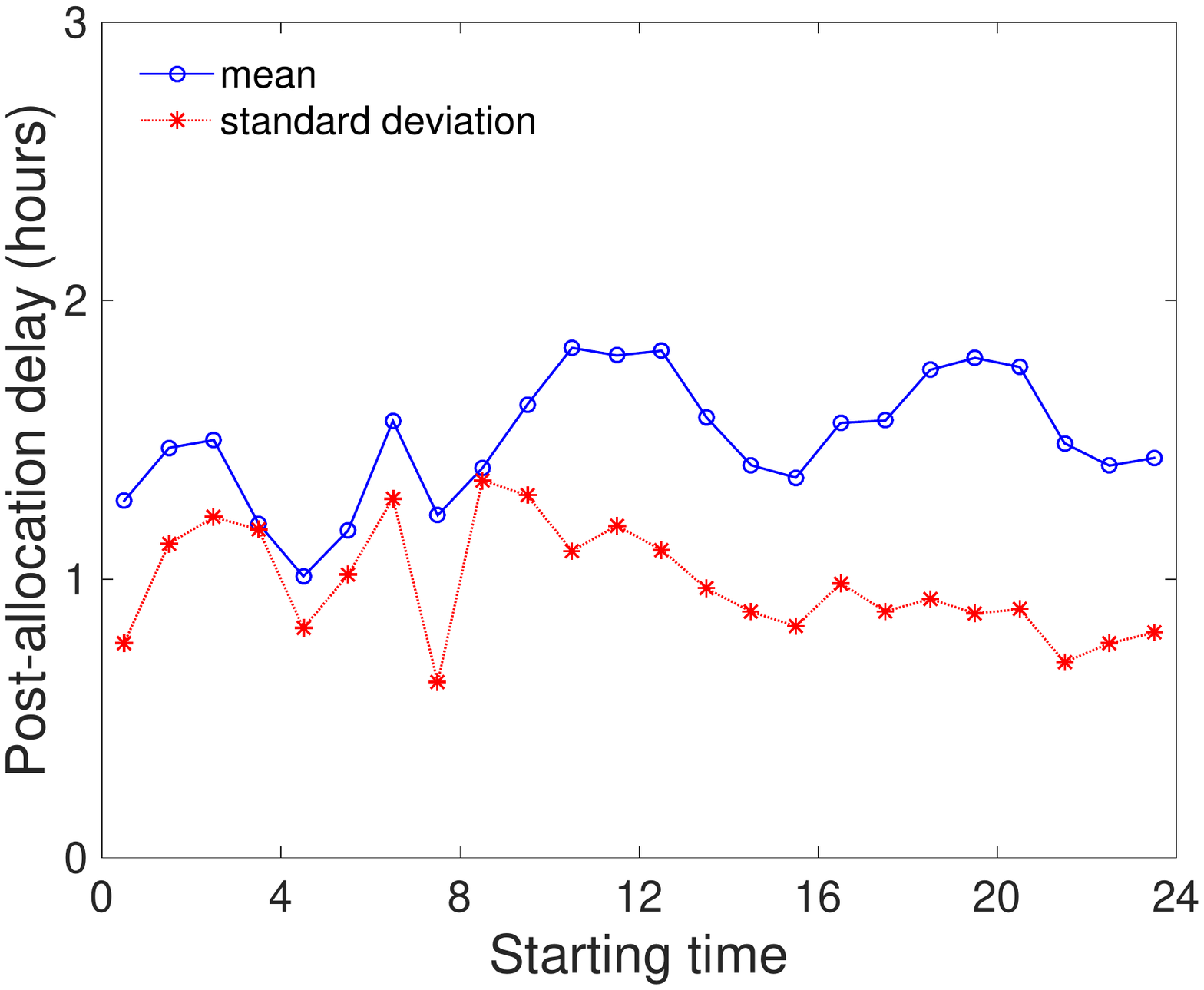}
\caption{Post-allocation delays}
\label{fig:post-allocation}
\end{subfigure}
\medskip
\caption{The means and standard deviations of pre- and post-allocation delays.}
\label{fig:allocation}
\end{figure}

We plot the means and standard deviations of pre- and post-allocation delays across a day obtained from the patient flow data in Figure~\ref{fig:allocation}. In the simulation model, we assume that these delays follow log-normal distributions with time-varying means and standard deviations specified in this figure.

In Table~\ref{tab:validate-delays}, we compare the mean waiting times from the ED to various bed pools obtained by simulation with those from the patient flow records. The simulation results generally agree with the results from the data for large bed pools (e.g., the pools of B2 and C beds in medicine and surgery). However, it is difficult for our simulation model to replicate mean waiting times for small bed pools. This is because the sample sizes of data about such pools are also small.

\begin{table}[t]
\centering
\caption{Mean lengths of stay (in days) in the six specialties from various sources, with corresponding standard deviations in the parentheses. A superscript $ \boldsymbol{*} $ denotes a sample of less than 5 admissions.}
\label{tab:LOS}
\medskip
\begin{tabular}{l l l l l l}
\toprule
&  \multicolumn{1}{c}{ED(am)} &  \multicolumn{1}{c}{ED(pm)} & \multicolumn{1}{c}{SDA} & \multicolumn{1}{c}{EL} & \multicolumn{1}{c}{SOC} \\ \midrule
Cardio. &  2.56 (3.23) & 3.10 (3.83) & 1.23 (1.16) & 2.21 (1.62)& 3.68 (4.79)  \\
Gastro. &  4.55 (4.37) & 4.70 (4.37) & NA & NA & NA \\
Geria. &  8.40 (9.29) & 8.90 (10.02) & NA & 19.30 (11.93)* & 11.86 (15.28)\\
Med. & 5.74 (6.87) & 6.16 (7.37)& 2.60 (2.60)*& 2.50 (2.28) & 5.74 (5.84)\\
Respir. & 3.87 (4.04) & 4.49 (4.68) & 1.33 (1.53)* & 1.33 (0.77)& 4.30 (3.02)\\
Surg. & 4.45 (5.44) & 5.03 (5.90) & 2.11 (2.42) & 2.71 (3.09) & 6.12 (7.07)\\
\bottomrule
\end{tabular}
\end{table}

\begin{table}[t]
\centering
\caption{The percentages of requested accommodation classes in the six specialties.}
\label{tab:accommodation}
\medskip
\begin{tabular} {l c c c c c c}
\toprule
& Cardio. & Gastro. & Geria. & Med. & Respir. & Surg. \\ \midrule
A1 &  \phantom{0}6\%   & \phantom{0}7\% & \phantom{0}4\% &  \phantom{0}5\%  & \phantom{0}4\% & \phantom{0}9\%   \\
B1 & \phantom{0}8\% & \phantom{0}9\% & \phantom{0}3\% & \phantom{0}7\% & \phantom{0}8\% & \phantom{0}8\% \\ 
B2 & 36\% & 39\% & 32\% & 32\% & 36\% & 39\%  \\
C & 50\% & 45\% & 61\% & 56\% & 52\% & 44\% \\ 
\bottomrule
\end{tabular}
\end{table}

\begin{table}[t]
\centering
\caption{Gender distributions in the six specialties.}
\label{tab:gender}
\medskip
\begin{tabular} {l c c c c c c}
\toprule
& Cardio. & Gastro. & Geria. & Med. & Respir. & Surg. \\ \midrule
Male &  64\% & 55\% & 34\% & 51\% & 63\% & 60\%  \\
Female & 36\% & 45\% & 66\% & 49\% & 37\% & 40\%  \\ 
\bottomrule
\end{tabular}
\end{table} 

\begin{table}[b]
\centering
\caption{The indexes of 34 bed pools in the hospital, with corresponding bed numbers in the parentheses.}
\label{tab:bed-pools}
\medskip
\begin{tabular}{l l l l l l l l l}
\toprule
& \multicolumn{1}{c}{Cardio.} & \multicolumn{1}{c}{Gastro.} & \multicolumn{1}{c}{Geria.} & \multicolumn{1}{c}{Respir.} & \multicolumn{1}{c}{Med.} & \multicolumn{1}{c}{Surg.} & \multicolumn{1}{c}{Surg. (EL)} & \multicolumn{1}{c}{Mixed} \\ \midrule
B1 Male & 1 (4) & NA & NA & NA & 9 (12) & 12 (18) & NA  & NA \\
B2 Male & 2 (24) & 3 (6) & 5 (4)& 7 (10) & 10 (26) & 13 (36) & 15 (14) & NA \\
C Male  & NA & 4 (8) & 6 (4) & 8 (12) & 11 (47) & 14 (42) & NA & 16 (19)\\
B1 Female  & 17 (4)  & NA & NA & NA & 25 (8) & 28 (16) & NA & NA \\
B2 Female & 18 (12) & 19 (6) & 21 (8) & 23 (12) & 26 (24) & 29 (20) & 31 (20) & NA \\
C Female & NA & 20 (8) & 22 (8) & 24 (9) & 27 (38) & 30 (36) & NA & 32 (16) \\
A1 & NA & NA & NA & NA & 33 (20) & 34 (20) & NA & NA \\
\bottomrule
\end{tabular}
\end{table}

\begin{table}[t]
\centering
\caption{Twenty-five types of male patients and their respective primary, preferred, and secondary bed pools.}
\label{tab:overflow-male}
\medskip
\begin{tabular}{l l l l l l}
\toprule
 & \multicolumn{1}{c}{Primary} & & \multicolumn{1}{c}{Preferred}  & & \multicolumn{1}{c}{Secondary} \\\midrule
Cardio. A1 & 33 & & NA & & 1, 9, 34   \\
Cardio. B1 & 1 & & NA & & 9, 33, 34 \\
Cardio. B2 & 2 & & NA & & 3, 4, 7, 8, 10, 11, 13, 16 \\
Cardio. C & 2 & & NA & & 3, 4, 7, 8, 10, 11, 14, 16\\
Gastro. A1 & 33 & & 9 & & 1, 12\\
Gastro. B1 & 9 & & 33 & & 12 \\
Gastro. B2 & 3, 7, 10 & & 4, 8, 11 & & 13, 16 \\
Gastro. C & 4, 8, 11 & &   3, 7, 10 & & 14, 16 \\
Geria. A1 & 33 & & NA  & & 9, 34, 12  \\
Geria. B1 & 9 & & NA & &  1, 12, 33 \\
Geria. B2 & 5 & & 6 & & 2, 3, 4 ,7, 8, 10, 11, 16  \\
Geria. C & 6 & & 5 & &  2, 3, 4, 7, 8, 10, 11, 16 \\
Respir. A1 & 33 & & 9 & & 34 \\
Respir. B1 & 9 & & 33 & & 1, 12\\
Respir. B2 & 3, 7, 10 & & 4, 8, 11 & & 2, 16 \\
Respir. C & 4, 8, 11 & & 3, 7, 10 & & 14, 16 \\
Med. A1 & 33 & & 9, 12 & & 34 \\
Med. B1 & 9 & & 33 & & 1, 12 \\
Med. B2 & 3, 7, 10 & &  4, 8, 11 & & 13, 14, 16\\
Med. C & 4, 8, 11 & & 3, 7, 10 & & 13, 14, 16 \\
Surg. A1 & 34 & & 12 & & 9 \\
Surg. B1 & 12 & & 34 & & 9 \\
Surg. B2 & 13 & & 14, 15 & & 10, 11, 16\\
Surg. C & 14 & & 13, 15 & & 10, 11, 16\\
Surg. B2 (EL) & 15 & & 14 & & 10, 11, 16\\
\bottomrule
\end{tabular}
\end{table}

\begin{table}[t]
\centering
\caption{Twenty-five types of female patients and their respective primary, preferred, and secondary bed pools.}
\label{tab:overflow-female}
\medskip
\begin{tabular}{l l l l l l}
\toprule
& \multicolumn{1}{c}{Primary} & & \multicolumn{1}{c}{Preferred}  & & \multicolumn{1}{c}{Secondary} \\\midrule
Cardio. A1 & 33 & & NA  & & 17, 25, 34 \\
Cardio. B1 & 17 & & NA  & & 25, 33, 34 \\
Cardio. B2 & 18 & & NA & & 19, 20, 23, 24, 26, 27, 29, 32 \\
Cardio. C & 18  & & NA & & 19, 20, 23, 24, 26, 27, 30, 32 \\
Gastro. A1 & 33 & & 25 & & 17, 28\\
Gastro. B1 & 25 & & 33 & & 28 \\
Gastro. B2 & 19, 23, 26 & & 20, 24, 27 & & 29, 32 \\
Gastro. C & 20, 24, 27 & & 19, 23, 26 & & 30, 32 \\
Geria. A1 & 33 & & NA  & & 25, 34, 28  \\
Geria. B1 & 25 & & NA & & 17, 28, 33 \\
Geria. B2 & 21 & & 22 & & 18, 19, 20, 23, 24, 26, 27, 32 \\
Geria. C & 22 & & 21 & & 18, 19, 20, 23, 24, 26, 27, 32 \\
Respir. A1 & 33 & & 25 & & 34 \\
Respir. B1 & 25 & & 33 & & 17, 28 \\
Respir. B2 & 19, 23, 26 & & 20, 24, 27 & & 18, 32 \\
Respir. C & 20, 24, 27 & & 19, 23, 26 & & 30, 32 \\
Med. A1 & 33 & & 25, 28 & & 34 \\
Med. B1 & 25 & & 33 & & 17, 28 \\
Med. B2 & 19, 23, 26 & &  20, 24, 27 & & 29, 30, 32 \\
Med. C & 20, 24, 27 & & 19, 23, 26 & & 29, 30, 32 \\
Surg. A1 & 34 & & 28 & & 25 \\
Surg. B1 & 28 & & 34 & & 25 \\
Surg. B2 & 29 & & 30, 31 & & 26, 27, 32\\
Surg. C & 30 & & 29, 31 & & 26, 27, 32\\
Surg. B2 (EL) & 31 & & 30 & & 26, 27, 32\\
\bottomrule
\end{tabular}
\end{table}

\begin{table}[t]
\centering
\caption{Estimates of ED patients' waiting times of various types, with 95\% confidence intervals.}
\label{tab:validate-delays}
\medskip
\begin{tabular}{l cc c cc}
\toprule
& \multicolumn{2}{c}{Male} & & \multicolumn{2}{c}{Female}\\ 
\cmidrule{2-3} \cmidrule{5-6}
& Simulation & Data &  & Simulation & Data \\ \midrule
Cardio. A1 & 3.78 $\pm$ 0.54 & 5.54 $\pm$ 1.91 & & 4.59 $\pm$ 0.74 & 6.08 $\pm$ 2.69 \\
Cardio. B1 & 3.48 $\pm$ 0.51 & 5.26 $\pm$ 2.00 & & 2.96 $\pm$ 0.69 & 4.65 $\pm$ 2.29 \\
Cardio. B2 & 6.03 $\pm$ 0.21 & 5.63 $\pm$ 0.74 & & 6.29 $\pm$ 0.29 & 6.43 $\pm$ 0.94 \\
Cardio. C & 6.13 $\pm$ 0.19 & 6.51 $\pm$ 0.68  & & 5.92 $\pm$ 0.26 & 6.89 $\pm$ 0.97 \\
Gastro. A1 & 2.65 $\pm$ 0.44 & 3.71 $\pm$ 2.66 & & 3.90 $\pm$ 0.49 & 5.84 $\pm$ 2.26 \\
Gastro. B1  & 3.42 $\pm$ 0.42 & 2.76 $\pm$ 2.31 & & 5.96 $\pm$ 0.47 & 4.60 $\pm$ 2.31 \\
Gastro. B2 & 7.91 $\pm$ 0.19 & 7.78 $\pm$ 1.05 & & 6.51 $\pm$ 0.21 & 7.09 $\pm$ 1.08 \\
Gastro. C & 5.72 $\pm$ 0.18 & 6.52 $\pm$ 0.91 & & 6.79 $\pm$ 0.19 & 6.71 $\pm$ 1.14 \\
Geria. A1 & 4.23 $\pm$ 1.07 & 7.36 $\pm$ 3.05 & & 4.05 $\pm$ 0.79 & 3.95 $\pm$ 3.56 \\
Geria. B1 & 3.33 $\pm$ 1.21 & 4.03 $\pm$ 3.94 & & 5.15 $\pm$ 0.81 & 6.73 $\pm$ 3.12 \\
Geria. B2 & 6.68 $\pm$ 0.38 & 7.42 $\pm$ 1.47 & & 7.51 $\pm$ 0.27 & 9.56 $\pm$ 1.00\\
Geria. C & 8.59 $\pm$ 0.27 & 8.12 $\pm$ 1.08 & & 9.38 $\pm$ 0.19 & 8.72 $\pm$ 0.75 \\
Respir. A1 & 4.11 $\pm$ 0.68 & 6.13 $\pm$ 2.79 & & 3.97 $\pm$ 0.92 & 5.34 $\pm$ 2.86 \\
Respir. B1 & 4.22 $\pm$ 0.49 & 4.70 $\pm$ 1.83 & & 5.54 $\pm$ 0.63 & 6.02 $\pm$ 1.88 \\
Respir. B2 & 8.80 $\pm$ 0.23 & 8.48 $\pm$ 0.81 & & 4.44 $\pm$ 0.31 & 6.99 $\pm$ 0.95 \\
Respir. C & 8.82 $\pm$ 0.19 & 7.17 $\pm$ 0.60 & & 6.48 $\pm$ 0.25 & 6.68 $\pm$ 0.89 \\	
Med. A1 & 4.09 $\pm$ 0.41 & 5.02 $\pm$ 1.18 & & 4.21 $\pm$ 0.43 & 5.85 $\pm$ 1.14 \\
Med. B1 & 4.05 $\pm$ 0.33 & 5.73 $\pm$ 1.01 & & 5.30 $\pm$ 0.33 & 4.90 $\pm$ 0.87\\
Med. B2 & 7.49 $\pm$ 0.15 & 7.59 $\pm$ 0.46 & & 7.64 $\pm$ 0.15 & 7.97 $\pm$ 0.43\\
Med. C & 5.76 $\pm$ 0.12 & 6.54 $\pm$ 0.33 & & 8.42 $\pm$ 0.12 & 8.15 $\pm$ 0.35 \\
Surg. A1 & 4.84 $\pm$ 0.30 & 5.33 $\pm$ 0.74 & & 5.10 $\pm$ 0.37 & 4.46 $\pm$ 1.06  \\
Surg. B1  & 2.99 $\pm$ 0.33 & 5.26 $\pm$ 0.90 & & 2.90 $\pm$ 0.40 & 5.43 $\pm$ 0.91 \\
Surg. B2  & 6.10 $\pm$ 0.15 & 6.58 $\pm$ 0.41 & & 8.72 $\pm$ 0.18 & 8.23 $\pm$ 0.45  \\
Surg. C  & 8.55 $\pm$ 0.13 & 7.75 $\pm$ 0.36 & & 6.23 $\pm$ 0.16 & 6.19 $\pm$ 0.46 \\	
\bottomrule
\end{tabular}
\end{table}

\end{document}